\documentclass[11pt,twoside]{article}

\usepackage{float}
\usepackage{tikz}
\usepackage[all]{xy}
\usepackage{wrapfig}
\usepackage{enumitem}
\usepackage{multirow, array} 
\usepackage{setspace}  

\usepackage{latexsym}
\usepackage{amssymb,amsbsy,amsmath,amsfonts,amssymb,amscd}
\usepackage{epsfig, graphicx}
\usepackage{hyperref}
\newcommand{\commentout}[1]{}
\newcommand{\R}{\mathbb{R}}
\newcommand{\N}{\mathbb{N}}

\newcommand {\Chi} {{\bf \raise 2pt \hbox{$\chi$}} }

\newcommand {\proof} {\noindent {\bf Proof}. }
\newcommand{\beq}{\begin{equation}}
\newcommand{\eeq}{\end{equation}}
\newcommand{\bea} {\begin{array}{rl}}
\newcommand{\eea} {\end{array}}
\newcommand{\bepa}{\left\{ \begin{array}{l}}
\newcommand{\eepa} {\end{array}\right.}
\newtheorem{theorem}{Theorem}[section]

\newtheorem{definition}[theorem]{Definition}
\newtheorem{remark}[theorem]{Remark}

\numberwithin{equation}{section}

\newcommand{\qed}{{ \hfill
                      {\unskip\kern 6pt\penalty 500 \raise -2pt\hbox{\vrule\vbox to 6pt{\hrule width 6pt
                      \vfill\hrule}\vrule} \par}   }}

\title{\Large \bf Blow-up, steady states 
and long time behaviour of
excitatory-inhibitory nonlinear neuron models}

\author{Mar\' {\i}a J. C\'aceres$^*$
\and Ricarda Schneider$^\dagger$
}

\date{}

\begin{document}
\doublespacing
\maketitle
\pagestyle{plain}
\pagenumbering{arabic}


\begin{abstract}
Excitatory and inhibitory nonlinear noisy leaky integrate and fire models
are often used to describe neural networks.
Recently, new mathematical results have provided a better understanding of 
them. 
It has been proved  that a fully excitatory network can
blow-up in finite time, while a fully inhibitory network 
has a global in time solution for any initial data.
A general description of the steady states 
of a purely excitatory or inhibitory network has been also given.
We extend this study to the system composed of an excitatory population
and an inhibitory one.  We prove that this system can also blow-up in 
finite time and analyse its steady states and long time behaviour.
Besides, we
illustrate our analytical description with some numerical results.
The main tools used to reach our aims are: the control of an exponential moment
for the blow-up results, a more complicate strategy than that considered
in \cite{CCP} for studying the number of steady states, entropy methods
combined with Poincar\'e inequalities for the long time behaviour and,
finally, high order numerical schemes together with parallel computation
techniques in order to obtain our numerical results.
\end{abstract}

 \rule{65mm}{0.1mm}

{\footnotesize
\hspace{0.5cm} \noindent \emph {2010 Mathematics Subject Classification.}  
35K60, 35Q92, 82C31, 82C32, 92B20

\hspace{0.5cm}\noindent \emph{Key words and phrases.} 
 Neural networks, Leaky integrate and fire models, noise, blow-up,
 steady states,

\hspace{0.5cm}\noindent entropy, long time behaviour.

\hspace{0.5cm}\noindent $^*$Departamento de Matem\'atica Aplicada,
Universidad de Granada, 18071 Granada, Spain.

 \hspace{0.5cm}\noindent $^\dagger$Departament de Matem\`atiques,
  Universitat Aut\`onoma de Barcelona,  E-08193 - Bellaterra, Spain 

\hspace{0.5cm}\noindent and Departamento de Matem\'atica Aplicada,
Universidad de Granada, 18071 Granada, Spain.

$\ $
}


\newpage

\section{Introduction}

One of the simplest self-contained mean field models for 
neural networks is the Network of Noisy Leaky Integrate and Fire 
(NNLIF) model. It is a mesoscopic description of the
 neuronal network, which gives,
at time $t$, the probability of finding a neuron with membrane potential 
$v$. There is a great deal of works about Integrate and Fire neuron 
models, 
\cite{BrHa,RBW,Touboul_AQIF, BrGe,RB,brunel,BrHa,Touboul_2008, Henry:13,
GK, T, G}. However, from a mathematical point of view the 
properties of the model are not well known yet. 
In the last years, some works have pointed to this direction;
in \cite{CGGS} some existence results have been proven: For a fully
inhibitory network there is global in time existence, while
for a purely excitatory network there is a global in time solution 
only if the firing rate is finite
for every time. These results are consistent with the fact that
in the excitatory case, solutions can blow-up in finite time 
if the  
value of the
 connectivity 
parameter
of the network is large enough or if
 the initial datum is concentrated close enough to the threshold potential,
 \cite{CCP,CP,Henry:13,delarue2015global,delarue2015particle}. 
Moreover, in \cite{CCP,CP}
the set of steady states for non coupled excitatory or inhibitory networks 
has been described, and the long time behaviour in the 
linear case has been studied. Later, in \cite{carrillo2014qualitative},
 for the nonlinear case with small connectivity parameters it has been proved 
exponential convergence to the unique steady state. 

In the present work we extend these results to the excitatory-inhibitory coupled NNLIF model. This model was also studied in \cite{brunel}, where time delay and refractory period 
were included. Here we focus on other aspects:  
We prove that, although in a purely inhibitor network the solutions are global in time (see \cite{CGGS}), in the presence
 of excitatory neurons 
the system can blow-up  in finite time. 
We also analyze the set of stationary states, which is more complicate than in the case of purely 
excitatory or inhibitory networks, and  prove exponential convergence to the unique steady
state when all the connectivity parameters are small.  
This exponencial convergence can be demonstrated,
in terms of the entropy method, since for this case the system
is a "small" perturbation of the linear one.
Finally,  the complexity of the coupled excitatory-inhibitory network is  numerically described.

The study developed  in \cite{brunel}, together with the results
in the present paper, show that this simple model 
(excitatory-inhibitory coupled NNLIF) could describe phenomena 
well known in neurophy\-sio\-lo\-gy: Synchronous and asynchronous states. 
As in \cite{brunel}
we will call  \emph{asynchronous} the states in which 
the firing rate tends to be constant in time and  
 \emph{synchronous} every other state.
Experimental and computational results exhibiting such phenomena can be 
found in \cite{brunel} and  references therein. Thus,
on one hand, when there are  asymptotically stable steady states 
there are asyncronous states, since the firing rate tends 
asymptotically to be constant in time.
Moreover, the presence of several steady states could provide
a rich behaviour of the network, since multi-stability phenomena 
could appear.
On the other hand, when the model does not have stable steady states, 
there are syncronous states. In this sense, the blow-up phenomenon 
could
be understood as a synchronization of part of the network, because
the firing rate diverges for a finite time. 
Possibly, this entails that a part of the network 
synchronizes, and thus fires at the same time.

\

Other PDE models describing the spiking neurons, which are related
to the Fokker-Planck system considered in the present work, are
based on Fokker-Planck equations including conductance variables,
 \cite{Caikinetic,CCTa, perthame2013voltage} (and references therein)
and on the time elapsed models \cite{PPD,PPD2,pakdaman2014adaptation}.
In \cite{dumont2015noisy} the authors study the connection between
the Fokker-Planck and the elapsed models.
From the microscopic point of view, 
the first family of PDE models is obtained assuming that 
the  spike trains follow a Poisson process, while
recently the elapsed models have been derived as mean-field
limits of Hawkes processes 
\cite{chevallier2015microscopic,chevallier2015mean}.

\

The structure of the paper is as follows: In Section \ref{sec:model},
the mathematical model and the notion of solution
considered are introduced. The finite time blow-up phenomenon is studied in 
Section \ref{sec:bu}, and the set of steady states as well as the long time
behaviour of the system are analyzed in Section \ref{sec:steadystates}.
Finally, in Section \ref{sec:numericalresults}  the analytical
results are illustrated numerically.


\section{The model and the definition of solution}
\label{sec:model}
We consider a neural network with $n$ neurons ($n_E$ excitatory
and $n_I$ inhibitory) described by 
the Integrate and Fire model, which depicts the activity of the 
membrane potential. 
The time evolution of the membrane potential $V_\alpha(t)$ of an 
inhibitory neuron ($\alpha=I$) or an excitatory one ($\alpha=E$) is given 
 by the following equation (see \cite{brunel, BrHa} for details)
\begin{eqnarray}
C_m \frac{dV^\alpha}{dt}(t)=-g_L(V^\alpha (t)-V_L)+I^\alpha(t), \label{ec1}
\end{eqnarray}
where $C_m$ is the capacitance of the membrane, $g_L$ is the
leak conductance, $V_L$ is the leak reversal potential
 and $I^\alpha(t)$ is the incoming synaptic current, 
which models all the interactions of the neuron with other neurons.
In the absence of  interactions with other neurons ($I^\alpha(t)=0$), 
 the  membrane potential  relaxes towards a resting value 
$V_L$. However, the interaction  with other neurons provokes the
neuron to fire, that is, it emits an action potential (spike) when 
$V^\alpha(t)$  reaches  its threshold or firing value $V_F$, and 
the membrane potential
relaxes to a reset value $V_R$. (Let us remark that $V_L<V_R<V_F$).
Each neuron receives 
$C_{ext}$ connections from excitatory neurons outside the network, 
and $C= C_E+C_I$ connections from neurons in the network;
$ C_E=\epsilon \, n_E$ from excitatory neurons 
and $C_I=\epsilon \, n_I$  from inhibitory neurons. 
 It is also supposed that $C_{ext}= C_E$.
These connections are assumed to be randomly choosen,
and the network  to be sparsely connected, namely, 
 $\epsilon=\frac{ C_E}{n_E}= \frac{C_I}{n_I}<< 1$, see \cite{brunel}. 
The synaptic current $I^\alpha(t)$ takes the form of the following
stochastic process 
\begin{eqnarray*}
I^\alpha(t)=J_E^\alpha \sum_{i=1}^{{\bar C}_E}\sum_j \delta (t-t_{Ej}^i)-
J_I^\alpha \sum_{i=1}^{C_I}\sum_j \delta (t-t_{Ij}^i), \qquad \alpha=E,I,
\end{eqnarray*}
where $t^i_{Ej}$ and $t^i_{Ij}$ are the times of the
$j^{th}$-spike coming from the $i^{th}$-presynaptic neuron for excitatory 
and inhibitory neurons, respectively, ${\bar C}_E= C_E+C_{ext}$,
and $J_k^\alpha$, 
for $\alpha,k=E,I$ are the strengths of the synapses.
The stochastic character is enclosed in the distribution of the spike 
times of the neurons. 
The spike trains of all neurons in the network are supposed to
be described by Poisson processes with a common instantaneous firing rate,
$\nu_\alpha(t)$, $\alpha=E,I$. These processes are supposed to be
 independent   \cite{brunel, CCP}. 
By using these hypotheses, the mean value of the current, 
$\mu
_C^{\alpha}(t)
$, 
and its variance, $\sigma_C^{\alpha 2}
(t)
$, take the form
\begin{eqnarray}
\mu^{\alpha}_C 
(t)
 &=& {C}_E J_E^\alpha \nu_E(t)-C_I J_I^\alpha \nu_I(t), \label{mu} \\
\sigma_C^{\alpha 2} 
(t)
 & =&{C}_E(J_E^\alpha)^2  \nu_E(t)
+
C_I (J_I^\alpha)^2 
\nu_I(t) \label{sigma}.
\end{eqnarray}
Many authors \cite{brunel, BrHa, mg,omurtag} 
then
 approximate 
the incoming synaptic current by a continuous in time stochastic 
process of Ornstein-Uhlenbeck type which has the same mean and variance 
as the Poissonian spike-train process. 
Specifically,  $I^\alpha(t)$ is approached by
\begin{eqnarray*}
I^\alpha(t)dt  \approx  \mu^\alpha_C 
(t)
\ dt+ \sigma_C^\alpha 
(t) \  d
B_t, \qquad \alpha=E,I,
\end{eqnarray*}
where $B_t$ is the standard Brownian motion.

Summing up, the approximation to the stochastic diferential equation 
model (\ref{ec1}), taking the voltage and time units so that $C_m=g_L=1$, 
finally yields
\begin{equation}
dV^\alpha  
(t)
 =  (-V^\alpha 
(t)
 + V_L+\mu_C^\alpha 
(t)
) \ dt+\sigma_C^\alpha 
(t)
 \ d B_t, 
\qquad V^\alpha\leq V_F, \qquad \alpha=E,I \label{ec2}
\end{equation}
with the jump process $V^\alpha(t_0^+)=V_R$, $V^\alpha(t_0^-)=V_F$, 
whenever at $t_0$ the voltage reaches the threshold value $V_F$.

The firing rate or
probability of firing per unit time of the Poissonian spike train, 
$\nu_\alpha
(t)
$, is calculated  in \cite{RBW} as
\begin{equation}
\nu_\alpha(t)= \nu_{\alpha,ext}+N_\alpha(t),  \qquad \alpha=E,I,
\nonumber
\end{equation}
where $\nu_{\alpha,ext}$ is the frequency of the 
external input
 and $N_\alpha(t)$ is the mean firing rate of the 
population $\alpha$. 
Also $\nu_{I,ext}=0$ since
the external connections are with excitatory neurons.

Going back to (\ref{ec2}), a system of
coupled partial differential equations  
for the evolution of the proba\-bi\-lity densities $\rho_\alpha(v,t)$ can be 
written, where $\rho_\alpha(v,t)$ denotes the probability of finding 
a neuron in the population $\alpha$, with a voltage $v \in (-\infty, V_F]$ 
at a time $t\ge 0$.
In  \cite{brunel, BrHa, mg,omurtag, Risken} a heuristic argument using 
It\^o's rule gives a system of coupled Fokker-Planck or backward Kolmogorov 
equations with sources
\begin{eqnarray}
\left \{ 
\begin{array}{l}
\frac{\partial\rho_I}{\partial t}(v,t)+
\frac{\partial}{\partial v} [h^I(v,N_E(t),N_I(t))
\rho_I(v,t)]-a_I(N_E(t),N_I(t))\frac{\partial^2\rho_I}{\partial v^2}(v,t)  
= N_I(t)\delta(v-V_R), 
\\
\\
\frac{\partial\rho_E}{\partial t}(v,t)+
\frac{\partial}{\partial v} [h^E(v,N_E(t),N_I(t))
\rho_E(v,t)]-a_E(N_E(t),N_I(t))\frac{\partial^2\rho_E}{\partial v^2}(v,t) 
= N_E(t)\delta(v-V_R),
\label{modelo}
\end{array} \right.
\end{eqnarray}
with $h^\alpha (v,N_E(t),N_I(t))  =  -v+ V_L+ \mu_C^\alpha$ and 
$ a_\alpha(N_E(t),N_I(t))=\frac{\sigma_C^{\alpha 2}}{2}$. 
The right hand sides in \eqref{modelo} represent the fact that when 
neurons reach the threshold potential $V_F$, they emit a spike over 
the network and reset their membrane potential to the reset value $V_R$. 
The system (\ref{modelo}) is completed with Dirichlet 
boundary conditions and an initial datum
\begin{eqnarray}
  \rho_\alpha(-\infty,t)=0, \
 \rho_\alpha(V_F,t)=0, \ \rho_\alpha(v,0)=\rho_\alpha^0(v) \ge 0,
 \quad   \alpha=E, I. \label{ec4}
\end{eqnarray}
In order to simplify the notation,
we denote $ {d}^\alpha_k=C_k(J_k^\alpha)^2\geq0 $ and 
$b_k^\alpha=C_kJ_k^\alpha \ge 0 $ for  $k, \alpha=E,I$,
and the variable $v$ is translated with the factor $V_L+b_E^E\, \nu_{E,ext}$.
Let us remark that we keep the same notation for the other involved 
values ($V_R, \ V_F$) and also $v$ for the new variable.
With the new voltage variable 
 and using expressions \eqref{mu} and \eqref{sigma} 
for $\mu_C^\alpha(t)$ and $\sigma_C^\alpha(t)$
, the drift and 
diffusion coefficients 
become
\begin{eqnarray}
h^\alpha (v,N_E(t),N_I(t))  & = &  
-v+b_E^\alpha N_E(t)-b_I^\alpha N_I(t)+(b_E^\alpha-b_E^E) \nu_{E ,ext}
, 
\label{drift}
\\ 
a_\alpha(N_E(t),N_I(t)) & = & {d}_E^{\alpha} \nu_{E,ext} + 
{d}_E^\alpha N_E(t)  
+
{d}_I^\alpha  N_I(t),  \quad \alpha=E,I. 
 \label{diffusion}
\end{eqnarray}
 The coupling of the system \eqref{modelo} is hidden in these two terms,
 since the mean firing rates $N_\alpha$ obey to
\begin{equation}
N_\alpha (t)  =  -a_\alpha (N_E(t), N_I(t))
\frac{\partial \rho_\alpha }{\partial v}(V_F,t) \ge 0, \quad \alpha=E,I. 
\label{N}
\end{equation}
Moreover, \eqref{N} gives rise to the nonlinearity of the system 
\eqref{modelo}, since firing rates  are defined in terms of boundary 
conditions on distribution functions $\rho_\alpha$.
On the other hand, since $\rho_E$ and $\rho_I$ represent probability densities,
the total mass should be conserved:
\begin{eqnarray*}
\int_{-\infty}^{V_F}\rho_{\alpha}(v,t) \ dv = 
\int_{-\infty}^{V_F}\rho_{\alpha}^0(v) \ dv = 1 \quad \forall \ t\ge 0, 
\quad \alpha=E,\, I.
\end{eqnarray*}

\

Before introducing the definition of solution considered in this work, let us
summarize some notations.
For $1\leq p < \infty$, the space of functions such that $f^p$ is 
integrable in $\Omega$ is denoted by $L^p(\Omega)$, 
$L^\infty (\Omega)$ is the space of essentially bounded functions in $\Omega$, 
$C^\infty(\Omega)$ is the set of infinitely differentiable functions in 
$\Omega$ and $L^1_{loc,+}(\Omega)$ denotes the set of non-negative 
functions that are locally integrable  in $\Omega$.
\begin{definition}
A weak solution of (\ref{modelo})-(\ref{N}) is a quadruple of 
nonnegative functions $(\rho_E, \rho_I, N_E, N_I)$ with 
$\rho_\alpha \in   L^\infty(\mathbb{R}^+;L^1_+((-\infty, V_F)))$
and  
$N_\alpha \in \ L^1_{loc,+}(\mathbb{R}^+) \ \forall \ \alpha=E,I$, satisfying
\begin{eqnarray}
\int_{0}^{T} \int_{-\infty}^{V_F} \rho_\alpha (v,t) 
\left[-\frac{\partial \phi}{\partial t} -
\frac{\partial \phi}{\partial v} h^\alpha(v,N_E(t),N_I(t))-
a_\alpha(N_E(t),N_I(t))\frac{\partial^2 \phi}{\partial v^2}\right] dv \ dt 
& \ & \ \label{debil} \\ 
= 
\int_{0}^{T} N_\alpha(t)[\phi(V_R,t)-\phi(V_F,t)]dt+ 
\int_{-\infty}^{V_F} \rho_\alpha^0(v)\phi(v,0) dv - \int_{-\infty}^{V_F}
\rho_\alpha(v,T)\phi(v,T) dv, & \ &  \alpha=E,I,
\nonumber
\end{eqnarray}
for any test function $\phi(v,t) \in  C^\infty((-\infty, V_F]  \times [0,T])$ 
such that $\frac{\partial^2 \phi}{\partial v^2}, 
\ v\frac{\partial \phi}{\partial v} \in L^\infty((-\infty, V_F) \times (0, T))$.
\end{definition}
Additionally, if  test functions of the form $\psi(t)\phi(v)$ are considered,
 the formulation  \eqref{debil} is equivalent to say that for all 
$\phi(v) \ \in \ C^\infty((-\infty,V_F])$ such that 
$v\frac{\partial \phi}{\partial v} \ \in \ L^\infty((-\infty,V_F))$
\begin{eqnarray}
\frac{d}{dt}\int_{-\infty}^{V_F}\phi(v)\rho_\alpha(v,t) \ dv  & = & 
\int_{-\infty}^{V_F} \left[ 
\frac{\partial \phi}{\partial v} h^\alpha(v,N_E(t),N_I(t))+
a_\alpha(N_E(t),N_I(t))\frac{\partial^2\phi}{\partial v^2} \right] 
\rho_\alpha(v,t) \ dv \nonumber \\
 \ & \ & + N_\alpha(t)\left[\phi(V_R)-\phi(V_F)\right] 
\label{debil2}
\end{eqnarray}
holds in the distributional sense for $\alpha = E,I$.  
Checking that weak solutions conserve the mass of the initial data 
is a straightforward computation after choosing $\phi=1$ in (\ref{debil2}), 
\begin{eqnarray*}
\int_{-\infty}^{V_F}\rho_{\alpha}(v,t) \ dv = 
\int_{-\infty}^{V_F}\rho_{\alpha}^0(v) \ dv = 1 \ \ \forall \ t\ge 0.
\end{eqnarray*}


\section{Finite time blow-up}
\label{sec:bu}
In \cite{CCP} and \cite{CP} it  was proved that weak solutions can blow-up 
in finite time for a purely excitatory network, when neurons
are considered or not to remain at a refractory state for a time. 
However, for a purely inhibitory network it was shown in \cite{GG09}
that weak solutions are global in time.
The following theorem claims that a network with excitatory and
 inhibitory neurons can blow up in finite time. 

We remark that the theorem is formulated in a more general setting
of drift terms $h^E$ than that considered in \eqref{drift}.
The diffusion term \eqref{diffusion} 
of the excitatory equation
is assumed to not vanish at any time.
For the inhibitory firing rate we assume \eqref{NI}, which
is satisfied, for instance,
 if $N_I(t)$ is bounded for every time.
This hypothesis should not be a strong constraint, because
in \cite{CGGS} it was proved, in the case of only one population
(in average excitatory or inhibitory), that if
the firing rate is bounded for every time, then there exists a
global solution in time. It could be natural to think that
an analogous criterion should hold in a coupled network, although its
proof seems much more complicated and remains as an open problem.
\begin{theorem}
Assume that
 \begin{equation}
h^E(v,N_E,N_I)+v \ge b_E^EN_E-b_I^EN_I, 
\label{hip1}
\end{equation}
\begin{equation}
 a_E(N_E,N_I)\ge a_m>0,  \label{hip2}
\end {equation}
$\forall v \; \in (-\infty,V_F]$ and $\forall \ N_I,N_E\ge0$.
Assume also that there exists some $M>0$ such that
\begin{equation}
\label{NI}
\int_0^tN_I(s) \ ds \le M\, t, \quad \quad \forall \,  t\ge 0.
\end{equation}
Then,
a weak solution to the system (\ref{modelo})-(\ref{N}) 
cannot be global in time because one
of
 the following reasons:
\begin{itemize}
\item  $b_E^E>0$ is large enough, for $\rho_E^0$ fixed.
\item $\rho_E^0$ is 'concentrated enough' around $V_F$:
\begin{eqnarray}
\int_{-\infty}^{V_F} e^{\mu v} \rho_E^0(v) \ dv \ge 
\frac{e^{\mu V_F}}{b_E^E \mu},
 \qquad \mbox{ for a certain } \mu>0 \label{cicritica}
\end{eqnarray}
and for $b_E^E>0$ fixed.
\end{itemize}
\label{th_blowup}
\end{theorem}
\proof
Using \eqref{debil2},
considering 
$\mu = \max\left( \frac{b_I^EM+2V_F
}{a_m}, 
\frac{1}{b_E^E} \right) $ and the multiplier $\phi(v)=e^{\mu v}$,
a weak solution $(\rho_E(v,t), \rho_I(v,t), N_E(t), N_I(t))$ 
satisfies the following inequality
\begin{eqnarray}
\frac{d}{dt}   \int_{-\infty}^{V_F} \phi(v)\rho_E(v,t)\, 
dv 
&\ge & \mu\int_{-\infty}^{V_F}
\phi(v)\left[b_E^EN_E(t)-b_I^EN_I(t)
-v\right]\rho_E(v,t)\, dv \nonumber\\
& &  + \
\mu^2a_m\int_{-\infty}^{V_F}\phi(v)\rho_E(v,t)\, dv + 
N_E(t)\left[\phi(V_R)-\phi(V_F)\right]\nonumber\\
&\ge & 
 \mu \left[b_E^EN_E(t)-b_I^EN_I(t)
-V_F+\mu a_m\right] 
\int_{-\infty}^{V_F}\phi(v)\rho_E(v,t) \, dv
\nonumber \\
& & -N_E(t)\phi(V_F),\nonumber
\end{eqnarray}
where assumptions \eqref{hip1}-\eqref{hip2} and the fact that 
$v \in (-\infty, V_F)$ and $N_E(t)\phi(V_R)>0$ were used.
This inequality,
Gronwall's lemma
and the definition of $\mu$ provide
the following inequality for the exponential moment
$M_\mu (t):= \int_{-\infty}^{V_F} \phi(v) \rho_E(v,t) \, dv$:
\begin{eqnarray}
M_\mu (t)  & \ge &
e^{\mu \int_0^t f(s) \ ds}
\left[ 
M_\mu (0)
-  \phi(V_F) \, \int_0^t N_E(s)  \,
e^{-\mu \int_0^s
f(z)
\ dz} \ ds
\right],
\nonumber
\end{eqnarray}
where $f(s)=b_E^EN_E(s)-b_I^E N_I(s)+\mu\, a_m-V_F$.
Using the definition of $\mu$ and \eqref{NI}, we notice that
$$
-\phi(V_F) \, \int_0^t N_E(s)  \,
e^{-\mu \int_0^s
f(z) \ dz} \ ds \ge 
-\phi(V_F) \, \int_0^t N_E(s)  \,
e^{-\mu \int_0^s \left[ b_E^E N_E(z)  + \mu a_m -V_F -M b_I^E \right]
\ dz} \ ds
$$
and after some computations the right hand side of the previous
inequality can be bounded by $-\frac{\phi(V_F)}{\mu b_E^E}$.
Finally, the following inequality holds
\begin{eqnarray}
M_\mu (t)  & \ge &
e^{\mu \int_0^t f(s) \ ds}
\left[ 
M_\mu (0)
-  \frac{\phi(V_F)}{\mu b_E^E}
\right].
\nonumber
\end{eqnarray}
We observe that if the initial state satisfies
\begin{eqnarray}
b_E^E \, \mu \, M_\mu(0)>\phi(V_F), \label{condicion}
\end{eqnarray}
then, denoting $K=M_\mu (0)-  \frac{\phi(V_F)}{\mu b_E^E}>0$,
\begin{equation}
\int_{-\infty}^{V_F}\phi(v)\rho_E(v,t)\, dv =M_\mu(t)
\ge K e^{\mu \int_0^t f(s) \ ds}, \quad  \forall \, t \ge 0.
\label{d}
\end{equation}
On the other hand, using the again definition of $\mu$ and \eqref{NI},
we observe that
$$
\mu \int_0^t f(s) \ ds \ge \mu \left[b_E^E\int_0^tN_E(s) \ ds+
\left(\mu a_m-V_F-Mb_I^E\right)t \right]
\ge \mu \, V_F t.
$$
Thus, $e^{\mu \int_0^t f(s)\ ds}\ge e^{\mu V_F t}$
and consequently, considering \eqref{d}, we obtain
\begin{eqnarray}
\int_{-\infty}^{V_F}\phi(v)\rho_E(v,t)\, dv=M_\mu (t)  & \ge &
K\, e^{\mu V_F t}.
\nonumber
\end{eqnarray}

On the other hand,  since $\rho_E(v,t)$ is a probability density and $\mu > 0$,
for all $ t\ge 0$: $\int_{-\infty}^{V_F}\phi(v)\rho_E(v,t)\, dv \le e^{\mu V_F}$,
which leads to a contradiction if the weak solution is
assumed to be global in time. 
Therefore, to conclude the proof there only remains to show  
 inequality (\ref{condicion}) in the two cases of the
theorem.
\begin{enumerate}
\item For a fixed initial datum and $b_E^E$ large enough, 
$\mu$, $M_\mu(0)$ 
and $\phi(V_F)$
are fixed, thus (\ref{condicion}) holds.
\item For $b_E^E >0$ fixed, if the initial data satisfy 
(\ref{cicritica}) then condition (\ref{condicion}) holds immediately.
Now, there  only remains to show that 
such initial data exist.

For that purpose we can approximate an initial Dirac 
mass at $V_F$ by smooth probability densities, so that 
$\rho_E^0 \simeq \delta (v-V_F)$. This gives the following condition 
\begin{displaymath}
e^{\mu V_F} \ge \frac{e^{\mu V_F}}{b_E^E \mu}, 
\end{displaymath}
which is satisfied if $\mu>\frac{1}{b_E^E}$. So, with our initial choice 
of $\mu$ we can ensure that the set of initial data we are looking for
 is not empty.  \qed
\end{enumerate}

\begin{remark}
Hypothesis \eqref{NI} could be relaxed by
$\int_0^tN_I(s) \ ds \le M\, t+ C\int_0^tN_E(s) \ ds$ (where
$C>0$).
Moreover,  using a priori estimates (as done  in \cite{CCP}) 
it could be proved that
$\int_0^tN_I(s) \ ds \le M\, (1+t)+ C\int_0^tN_E(s) \ ds$, 
which seems not to be enough to reach the whole result.
Precisely, it yields the blow-up for fixed initial data and 
large $b_E^E$, but not
for fixed $b_E^E$ and  concentrated initial data.
\end{remark}

Theorem \ref{th_blowup} shows that 
 blow-up occurs when the connectivity
parameter of the excitatory to excitatory synapses, $b_E^E$, is large enough or 
if initially there are many excitatory neurons with a
voltage value of their membrane potential very close to the threshold value, $V_F$. 
From a biological perspective, and in view
of the numerical results in
Section \ref{sec:numericalresults} 
(Figs. \ref{blowup_bEE}, \ref{blowup_ci} and 
\ref{blowup_bEE_bIE}),
 in the first case, blow-up appears due
 to the strong influence of the excitatory population on the behavior of 
the network, poducing the incontrolled growth of the firing rate
at finite time. 
In the second case, even with  small  connectivity
parameter  $b_E^E$, the abundance of  excitatory neurons
with membrane potential voltage values sufficiently close to $V_F$ causes the firing rate 
to diverge in finite time. For a microscopic description, at the
level of individual neurons, we refer to \cite{delarue2015particle}
and \cite{delarue2015global} where the blow-up phenomenon is also analysed.

\

As mentioned 
above
it was proved  in \cite{CP} that 
fully excitatory networks  can blow-up in finite time, when one includes 
the refractory state.
This result can be extended to the case of excitatory-inhibitory 
networks by following the same ideas of Theorem \ref{th_blowup}.


\section{Steady states and long time behavior}
\subsection{Steady states}

\label{sec:steadystates}
For excitatory and inhibitory networks, 
the study of their steady states follows a similar strategy
 to that for the fully excitatory or inhibitory cases. 
However, for coupled networks the system solved by the stationary 
solutions is much more complicate, as we will see throughout this section.

As in \cite{CCP}, let us search for continuous stationary solutions 
$(\rho_E,\rho_I)$ of (\ref{modelo}) such that $\rho_E, \ \rho_I$ are 
$C^1$ regular, except possibly at $V=V_R$ where they are Lipschitz,
and satisfy the following identity:
\begin{equation}
\frac{\partial}{\partial v}[h^\alpha(v) \rho_\alpha(v)-
a_\alpha(N_E,N_I)\frac{\partial 
\rho_\alpha}{\partial v}(v)-N_\alpha H(v-V_R)]=0,
\quad \alpha=E,I,
\nonumber
\end{equation}
in the sense of distributions,
where $H$ denotes the Heaviside function.
Considering $h^\alpha(v, N_E,N_I)=V_0^\alpha(N_E,N_I)-v$, 
with
$V_0^\alpha (N_E,N_I)=b_E^\alpha N_E - b_I^\alpha N_I + 
(b_E^\alpha-b_E^E)\nu_{E,ext}$,
the definition of 
$N_\alpha$ in (\ref{N}) and the  Dirichlet boundary conditions in (\ref{ec4}),
we find the following initial value problem
for every $\alpha=E,I$:
\begin{equation}
\left \{ 
\begin{array}{l}
\frac{\partial \rho_\alpha}{\partial v}(v)  =  
\frac{V_0^\alpha(N_E,N_I)-v}{a_\alpha(N_E,N_I)}
\rho_\alpha(v)-\frac{N_\alpha H(v-V_R)}{a_\alpha(N_E,N_I)}, \\
\rho_\alpha(V_F)  =  0,
\nonumber
\end{array} \right.
\end{equation}
whose solutions are of the form
\begin{equation}
\rho_\alpha (v) = 
\frac{N_\alpha}{a_\alpha(N_E,N_I)}
e^{-\frac{(v-V_0^\alpha(N_E,N_I))^2}
{2a_\alpha(N_E,N_I)}}\int_{\max(v,V_R)}^{V_F}
e^{\frac{(w-V_0^\alpha(N_E,N_I))^2}{2a_\alpha(N_E,N_I)}}dw, 
\quad \alpha=E,I. \label{sol_eq}
\end{equation}
The expression  \eqref{sol_eq} is not an explicit formula
for $\rho_\alpha$, since the right hand side depends on $N_\alpha$, but
 provides a system of implicit equations for  $N_\alpha$
\begin{equation}
\frac{a_\alpha(N_E,N_I)}{N_\alpha}  = 
 \int_{-\infty}^{V_F}e^{-\frac{(v-V_0^\alpha(N_E,N_I))^2}{2a_\alpha(N_E,N_I)}} 
\left[ \int_{\max{(v,V_R)}}^{V_F}
e^{\frac{(w-V_0^\alpha(N_E,N_I))^2}{2a_\alpha(N_E,N_I)}}dw \right] dv, 
\qquad  \alpha = E,I, 
 \label{ecN}
\end{equation}
for which the conservation of mass 
($\int_{-\infty}^{V_F} \rho_\alpha(v,t)dv = 1$) has been used.
Therefore, the stationary solutions $(\rho_E,\rho_I)$ have the profile 
(\ref{sol_eq}), where $(N_E, N_I) $ are positive solutions of the implicit 
system (\ref{ecN}), and there are as many as solutions of (\ref{ecN}).
 Of course, in the particular case of a linear system,
 that is $V_0^\alpha(N_E,N_I)=0$ and $a_\alpha(N_E,N_I)$ independent
of the firing rates, there is a unique steady state.

Determining the number of solutions of the implicit system (\ref{ecN})
is a difficult task. In this section we find some conditions on the 
parameters of the model in order to reach this goal. 
Firstly, 
we consider the following change of variables and notation:
\begin{eqnarray}
z&\! \! \! \! \! =\! \! \! \! \! &\frac{v-V_0^E(N_E,N_I)}{\sqrt{a_E (N_E,N_I)}},  
\ u=\frac{w-V_0^E(N_E,N_I)}{\sqrt{a_E(N_E,N_I)}},  
\ w_F:=\frac{V_F-V_0^E(N_E,N_I)}{\sqrt{a_E(N_E,N_I)}},   \
w_R:=\frac{V_R-V_0^E(N_E,N_I)}{\sqrt{a_E(N_E,N_I)}},
\nonumber \\
\tilde{z}&\! \! \! \! \! =\! \! \! \! \! &\frac{v-V_0^I(N_E,N_I)}{\sqrt{a_I(N_E,N_I)}},  \
 \tilde{u}=\frac{w-V_0^I(N_E,N_I)}{\sqrt{a_I(N_E,N_I)}},  \
 \tilde{w}_F:=\frac{V_F-V_0^I(N_E,N_I)}{\sqrt{a_I(N_E,N_I)}},   \
 \tilde{w}_R:=\frac{V_R-V_0^I(N_E,N_I)}{\sqrt{a_I(N_E,N_I)}}.
\nonumber
\end{eqnarray}
With these new variables, the system (\ref{ecN}) is rewritten as
\begin{eqnarray}
\frac{1}{N_E} & = & I_1(N_E,N_I),
\ \mbox{where} \ I_1(N_E,N_I)  =  
\int_{-\infty}^{w_F}e^{-\frac{z^2}{2}}
\int_{\max(z,w_R)}^{w_F}e^{\frac{u^2}{2}}du \ dz,
 \nonumber \\
\frac{1}{N_I} & = & I_2(N_E,N_I), 
\ \mbox{where} \
I_2(N_E,N_I) = \int_{-\infty}^{\tilde{w}_F}
e^{-\frac{z^2}{2}}\int_{\max(z,\tilde{w}_R)}^{\tilde{w}_F}
e^{\frac{u^2}{2}}du \ dz.
\label{ecN2}
\end{eqnarray}
Now, using 
$s=\frac{z-u}{2}$ and $\tilde{s}=\frac{z+u}{2}$,
 the functions
$I_1$ and $I_2$ can be formulated as
\begin{eqnarray}
I_1(N_E,N_I)=
\int_0^\infty 
\frac{e^{-\frac{s^2}{2}}}{s}(e^{s \, w_F}-e^{s \, w_R}) \ ds,
 \nonumber\\
I_2(N_E,N_I)=\int_0^\infty \frac{e^{-\frac{s^2}{2}}}{s}
(e^{s \, \tilde{w}_F}-e^{s \, \tilde{w}_R}) \ ds. \nonumber
\end{eqnarray}
When $b_I^E=b_E^I=0$ the equations are uncoupled and the number of 
steady states can be determined in the same way as in \cite{CCP}, 
depending on the values of  $b_E^E$, since for the inhibitory equation 
there is always a unique steady state.
The following theorem establishes a classification of 
the number of steady states
in terms of the model parameters, in the case that the connectivity
parameters $b_I^E$ and $b_E^I$ do not vanish, and in comparison 
with the values of the pure connectivity parameters $b_E^E$ and $b_I^I$.

\begin{theorem}
Assume that the connectivity
parameters $b_I^E$ and $b_E^I$ do not vanish ($b_I^E, \, b_E^I>0$), 
$a_\alpha$ is independent of 
$N_E$ and $N_I$, $a_\alpha(N_E, N_I)=a_\alpha$,
 and $h^\alpha(v,N_E,N_I)=V_0^\alpha (N_E,N_I)-v$ with
$V_0^\alpha (N_E,N_I)=b_E^\alpha N_E - b_I^\alpha N_I + 
(b_E^\alpha-b_E^E) 
\nu_{E,ext}
$ for all $\alpha = E,I$.
Then:
\begin{enumerate}
\item There is an even number of steady states or there are no steady states 
for  
\eqref{modelo}-\eqref{N}
if
\begin{equation}
(V_F-V_R)^2< (V_F-V_R)(b_E^E-b_I^I) + 
 b_E^Eb_I^I-b_I^Eb_E^I.
\label{condicion21}
\end{equation}
If $b_E^E$ 
is large enough in comparison
with the rest of connectivity para\-me\-ters, there are no steady states.
Specifically,  there are no steady states 
 if both \eqref{condicion21}  and
\begin{equation}
 \max \left\{I_1(0)2\left(V_F+b_I^EN_I(\bar N_E)\right), 
\frac{b_E^Ib_I^E}{b_I^I}, 2(V_F-V_R)\right\}<b_E^E
\label{condicion3}
\end{equation} 
are fulfilled, where
\begin{equation}
\bar N_E =\max
\left\{N_E^*\ge 0:N_E^*=\frac{2(b_E^IN_I(N_E^*)+V_F)}{b_E^E}\right\}.
\label{eleccion_N}
\end{equation}
\item There is an odd number of steady states for  
\eqref{modelo}-\eqref{N} if
\begin{equation}
(V_F-V_R)(b_E^E-b_I^I) + b_E^Eb_I^I-b_I^Eb_E^I
<(V_F-V_R)^2.
\label{condicion21bis}
\end{equation}
If $b_E^E$  is small enough in comparison
with the rest of connectivity para\-me\-ters, there is a unique steady 
state.
\end{enumerate}
\label{th_eq}
\end{theorem}

\proof 
The proof reduces to study the existence of solutions to
 (\ref{ecN2}). It is organized in several steps. 
 
\emph{Step 1. To prove that  $\frac{1}{N_I}=I_2(N_E, N_I)$
admits a unique solution $N_I(N_E)$ for  $N_E$ fixed.}  

\begin{figure}[H]
\begin{center}
\includegraphics[height=6cm]{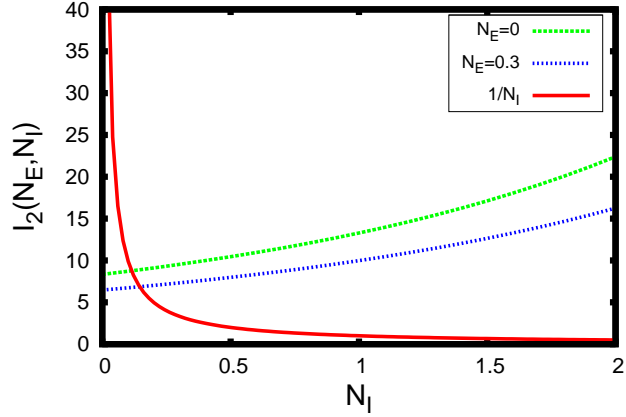}
\caption{The function 
$I_2$ in terms of $N_I$, for different values of $N_E$ fixed. 
$I_2$ is an increasing function on $N_I$ and 
decreasing on $N_E$.
}
\label{figI}
\end{center}
\end{figure}

Given  $N_E  \in [0, \infty)$ and following analogous ideas 
as in \cite{CCP}, it is easy to observe that there is a unique 
solution $N_I(N_E)$ to
\begin{eqnarray}
\frac{1}{N_I}=I_2(N_E, N_I) \label{I2}.
\end{eqnarray}
This fact is a consequence of the following properties of $I_2$:
\begin{enumerate}
\item  $I_2(N_E, N_I)$ is $C^\infty$ in both variables.
\item  For every $N_E$ fixed, $I_2(N_E, N_I)$ is an increasing strictly 
convex function on $N_I$, since for all integers $k\ge1$
\begin{equation}
\frac{\partial^k I_2}{\partial N_I^k} 
=
\left( \frac{b_I^I}{\sqrt{a_I}}\right)^k\int_0^\infty 
e^{-\frac{{s}^2}{2}}s^{k-1}(e^{s\tilde{w}_F}-e^{s\tilde{w}_R}) \ ds.
\label{parcial_NI}
\end{equation}
Thus, $\displaystyle \lim_{N_I \rightarrow \infty} I_2(N_E,N_I)=\infty$
for every $N_E$ fixed.
\item If we consider $N_I  \in  [0, \infty )$,   $I_2(N_E, N_I)$ is a
 decreasing convex function on $N_E$, since for all integers $k\ge1$
\begin{equation}
\frac{\partial^k I_2}{\partial N_E^k} =(
-1)^k\left( \frac{b_E^I}{\sqrt{a_I}}\right)^k
\int_0^\infty e^{-\frac{{s}^2}{2}}s^{k-1}
(e^{s\tilde{w}_F}
-e^{s\tilde{w}_R}) \ ds.
\label{parcial_NE}
\end{equation}
Thus, $\displaystyle \lim_{N_E \rightarrow \infty} I_2(N_E,N_I)=0$,
for every $N_I$ fixed.

\item Using expression \eqref{ecN2}  for
 $I_2$, we have $I_2(N_E, 0)<\infty$, for every $N_E$ fixed, 
since 
\begin{equation}
 I_2(N_E, 0)    =  
\int_{-\infty}^{\tilde{w}_F(0)}
e^{-\frac{z^2}{2}}
\int_{\max(z, \tilde{w}_R(0))}^{\tilde{w}_F(0)} 
e^{\frac{u^2}{2}}  \,  du \, dz
\le 
\sqrt{2 \pi} \left( \frac{V_F-V_R}{\sqrt{a_I}}\right)  
e^{\frac{m}{2a_I }},
\nonumber
\end{equation}
where
$m= \max 
\{(V_F-b_E^IN_E-(b_E^I-b_E^E)\nu_{E,ext})^2,
(V_R-b_E^IN_E-(b_E^I-b_E^E)\nu_{E,ext})^2 \}$.
\end{enumerate}
Figure \ref{figI} depicts the function 
$I_2$ in terms of $N_I$ for different values of $N_E$ fixed.
In this figure, the properties of $I_2$  enumerated above 
 can be observed.
\emph{Step 2. Properties of $N_I(N_E)$.} 

Obtaining an analytical expression of $N_I(N_E)$ is not
easy. However, some properties of $N_I(N_E)$ are enough to prove
the theorem:
\begin{enumerate}
\item $N_I(N_E)>0$, $\forall$ $N_E \in [0, \infty)$.
\item $N_I(N_E)$ is an increasing function, $\forall \ N_E \in [0,\infty)$,
since
\begin{equation}
N_I'(N_E)=
\frac{-N_I^2(N_E)\frac{\partial I_2}{\partial N_E}(N_E,N_I(N_E))}
{1+N_I^2(N_E)\frac{\partial I_2}{\partial N_I}(N_E,N_I(N_E))}
\label{N'}
\end{equation}
is nonnegative.
Expression \eqref{N'} is obtained by 
differentiating $N_I(N_E)I_2(N_E,N_I(N_E))=1$
with  respect to $N_E$.
\item $\displaystyle \lim_{N_E \rightarrow \infty}N_I(N_E)=\infty$,
because $N_I(N_E)$ is a positive increasing function,
thus its limit could be infinite or a constant $C>0$.
The  finite limit leads  to a contradiction
since $\frac{1}{N_I(N_E)}=I_2(N_E,N_I(N_E))$, and then
$\displaystyle \lim_{N_E \rightarrow \infty}\frac{1}{N_I(N_E)} =  \frac{1}{C}$,
while $\displaystyle \lim_{N_E \rightarrow \infty}I_2(N_E,N_I(N_E))=0$
because $\displaystyle \lim_{N_E\to \infty} I_2(N_E,C)=0$.
\item $0<N'_I(N_E)< \frac{b_E^I}{b_I^I}$.
It is a consequence of the fact that 
\begin{equation}
N_I'(N_E)=
\frac{b_E^I N_I^2(N_E) I(N_E)}
{\sqrt{a_I}+ b_I^I N_I^2(N_E)I(N_E)},
\label{N'bis}
\end{equation}
 where
$
I(N_E)=
\int_0^\infty e^{-s^2/2}
e^{\frac{-(b_E^IN_E-b_I^IN_I(N_E)+(b_E^I-b_E^E)\nu_{E,ext})s}{\sqrt{a_I}}}
\left( e^{sV_F/\sqrt{a_I}}-e^{sV_R/\sqrt{a_I}} \right) \ ds.
$

Expression \eqref{N'bis} is obtained using \eqref{parcial_NI}, \eqref{parcial_NE},  \eqref{N'} and the
definition of  $\tilde{w}_F$ and  $\tilde{w}_R$.
\item $\displaystyle \lim_{N_E\to \infty} N_I' (N_E)=\frac{b_E^IL}{\sqrt{a_I}+b_I^I\, L}$,
where
$\displaystyle\lim_{N_E \to \infty} N_I^2(N_E) I(N_E)=L$.
\item $L=\frac{\sqrt{a_I}}{V_F-V_R}$. This limit is obtained after some
tedious calculations, using a Taylor expansion of 
$e^{sV_F/\sqrt{a_I}}-e^{sV_R/\sqrt{a_I}}$ at $s=0$ and  computing 
the following limit:
$$
\lim_{N\to \infty} N^2
\int_0^\infty e^{-s^2/2} e^{-(b\, N- k)}s^n \, ds=
\left\{
\begin{array}{l}
\frac{1}{b^2} \quad \mbox{if} \quad n=1
\\
0 \quad \mbox{if} \quad n>1,
\end{array}
\right.
$$
for every $k\in \R$, $b>0$ and $n\in \N$.
\item As a consequence of the two previous properties, 
$\displaystyle \lim_{N_E \to \infty} N_I'(N_E)=\frac{b_E^I}{V_F-V_R+b_I^I}$.

The limit of $N_I'(N_E)$ ensures that the rate of
increase for $N_I(N_E)$ (in point 
3 of Step 2 it was proved that $N_I(N_E)$ tends to infinity) is 
 constant at infinity, and can be controlled in terms of 
the parameters $b_E^I$, $b_I^I$, $V_R$ and $V_F$.
\end{enumerate}

\emph{Step 3. Properties of the function 
$\mathcal{F}(N_E)=N_EI_1(N_E,N_I(N_E))$.}

For every fixed $N_E\in [0, \infty)$ it is obtained  $N_I(N_E)$, 
the unique solution
to  (\ref{I2}) (cf. Step 1). 
To conclude the proof there remains
   to determine the number of solutions to 
\begin{eqnarray}
N_EI_1(N_E, N_I(N_E))=1. \label{ecF}
\end{eqnarray}
With this aim, we define $\mathcal{F}(N_E)=N_EI_1(N_E,N_I(N_E))$. 
Depending on the properties of $\mathcal{F}$ we can find a different
number of solutions to (\ref{ecF}). First,
following analogous ideas as in \cite{CCP}, let us show
some properties of  $I_1$:

\begin{enumerate}

\item  $I_1(N_E, N_I)$ is $C^\infty$ in both variables.

\item For every $N_E \in [0,\infty)$ fixed,
 $I_1(N_E, N_I)$ is an increasing strictly convex 
function on $N_I$, since for all integers $k\ge1$
\begin{equation}
\frac{\partial^k I_1}{\partial N_I^k} =
\left( \frac{b_E^I}{\sqrt{a_E}}\right)^k
\int_0^\infty e^{-\frac{{s}^2}{2}}s^{k-1}(e^{s w_F}-e^{sw_R}) \ ds.
\nonumber
\end{equation}
Thus, $\displaystyle \lim_{N_I \to \infty} I_1(N_E,N_I)=\infty$,
for every $N_E \in [0,\infty)$ fixed.
\item For every $N_I \in [0, \infty )$ fixed, 
  $I_1(N_E, N_I)$ is a decreasing convex function on $N_E$, 
since for all integers $k\ge1$
\begin{equation}
\frac{\partial^k I_1}{\partial N_E^k} 
=(-1)^k\left( \frac{b_E^E}{\sqrt{a_E}}\right)^k
\int_0^\infty e^{-\frac{{s}^2}{2}}s^{k-1}(e^{sw_F}-e^{sw_R}) \ ds.
\nonumber
\end{equation}
Thus, $\displaystyle \lim_{N_E \to \infty} I_1(N_E,N_I)=0$,
 for every $N_I \in [0, \infty )$ fixed.

\item $I_1(N_E):=I_1(N_E, N_I(N_E))$ has the following pro\-per\-ties 
of monotonicity:
\begin{enumerate}[label*=\arabic*.]
\item If $b_E^I\, b_I^E <b_E^E\, b_I^I$, then $I_1(N_E)$ is a 
decreasing function with $ \displaystyle \lim_{N_E \to \infty} I_1(N_E)=0$.
\item If $b_E^I\, b_I^E >b_E^E\, b_I^I$  
\begin{enumerate}[label*=\arabic*.]
\item and 
$\frac{b_E^I}{V_F-V_R+b_I^I} <\frac{b_E^E}{b_I^E}$, 
then $I_1(N_E)$ decreases for $N_E$ large enough and
 $ \displaystyle \lim_{N_E \to \infty} I_1(N_E)=0$.
\item and $\frac{b_E^I}{V_F-V_R+b_I^I} >\frac{b_E^E}{b_I^E}$, then 
$I_1(N_E)$ increases for $N_E$ large enough and $ \displaystyle
 \lim_{N_E \to \infty} I_1(N_E)=\infty$.

If $b_E^E$ is small enough, such that $b_E^E<b_I^E\, N_I'(N_E)$ for
all $N_E\ge 0$,  $I_1(N_E)$ is an increasing function and 
$\displaystyle \lim_{N_E \to \infty} I_1(N_E)=\infty$.
\end{enumerate}
\end{enumerate}
These properties are proved using 
\begin{equation}
\frac{d}{d N_E} I_1(N_E, N_I(N_E))  =  
\left(
\frac{-b_E^E}{\sqrt{a_E}}
+\frac{b_I^E}{\sqrt{a_E}}N'_I(N_E)\right)
\int_0^\infty e^{\frac{-s^2}{2}}(e^{s w_F}-e^{s w_R})ds . 
\nonumber
\end{equation}
Therefore,  we have that
$I_1(N_E)=I_1(N_E, N_I(N_E))$ decreases 
iff $N'_I(N_E) < \frac{b_E^E}{b_I^E}$. Consequently,
 it
increases  iff $N'_I(N_E) > \frac{b_E^E}{b_I^E}$. 
Now, using properties 4. and 7. of 
$N_I(N_E)$  
the monotonicity properties of $I_1(N_E)$ are immediate.
\end{enumerate}
Taking into account the previous points,
the following properties of  $\mathcal{F}$ are obtained:
\begin{enumerate}
\item $\mathcal{F}(0)=0$.
\item Monotonicity. 
\begin{enumerate}[label*=\arabic*.]
\item If $b_E^Ib_I^E<b_E^Eb_I^I$ or if  $b_E^Ib_I^E>b_E^Eb_I^I$ and 
$\frac{b_E^I}{V_F-V_R+b_I^I}<\frac{b_E^E}{b_I^E}$, then 
 for $N_E$ large enough, $\mathcal{F}$ is a
 decreasing function.
\item If $b_E^Ib_I^E>b_E^Eb_I^I$ and 
$
\frac{b_E^I}{V_F-V_R+b_I^I} >\frac{b_E^E}{b_I^E}$, 
then  for $N_E$ large enough, $\mathcal{F}$ is an increasing function.
Note that if $b_E^E$ is small enough such that $b_E^E<b_I^E N_I'(N_E)$
for all $N_E \ge 0$, 
then $\mathcal{F}$ is an increasing function for all $N_E\ge 0$.
\item Close to $N_E=0$, $\mathcal{F}$ increases.
\end{enumerate}
The monotonicity of $\mathcal{F}$ is given by the sign of its derivative
\begin{equation}
\mathcal{F}'(N_E)= I_1(N_E)+ 
N_E\left[-\frac{b_E^E}{\sqrt{a_E}}+\frac{b_I^E}{\sqrt{a_E}}N_I'(N_E)\right]
\int_0^\infty e^{\frac{-s^2}{2}}(e^{sw_F}-e^{sw_R}) \ ds.
\nonumber
\end{equation}
We observe that
$
\lim_{N_E\to \infty}\mathcal{F}'(N_E)= 
\lim_{N_E\to \infty}\left(I_1(N_E)+ 
k\, N_E  \int_0^\infty e^{\frac{-s^2}{2}}(e^{sw_F}-e^{sw_R})ds
\right)$, where 
$k=\frac{(b_E^Ib_I^E-b_E^Eb_I^I)-b_E^E\,\left(V_F-V_R\right)}
{\sqrt{a_E}(b_I^I+V_F-V_R)}$. So, the monotonicity properties 
of $\mathcal{F}$ can be proven as follows.
\begin{enumerate}[label*=\arabic*.]
\item If $b_E^Ib_I^E<b_E^Eb_I^I$ or if  $b_E^Ib_I^E>b_E^Eb_I^I$ and 
$\frac{b_E^I}{V_F-V_R+b_I^I}<\frac{b_E^E}{b_I^E}$,  then
$\displaystyle \lim_{N_E\to \infty}\mathcal{F}'(N_E)=0_-$.
\item If $b_E^Ib_I^E>b_E^Eb_I^I$ and 
$\frac{b_E^I}{V_F-V_R+b_I^I}>\frac{b_E^E}{b_I^E}$, then
$\displaystyle \lim_{N_E\to \infty}\mathcal{F}'(N_E)=\infty$.
\item Since $\mathcal{F}(0)=0$ and $\mathcal{F}(N_E)>0$, $\mathcal{F}$ 
increases close to 0.
\end{enumerate}

\item Limit of $\mathcal{F}$. 
\begin{enumerate}[label*=\arabic*.]
\item If $b_E^Ib_I^E<b_E^Eb_I^I$ or if $b_E^Ib_I^E>b_E^Eb_I^I$ 
and $\frac{b_E^I}{V_F-V_R+b_I^I}
<\frac{b_E^E}{b_I^E}$, then
\begin{equation}
\lim_{N_E\to \infty}\mathcal{F}(N_E)= 
\frac{(V_F-V_R)(V_F-V_R+b_I^I)}{b_E^E(V_F-V_R)+b_E^Eb_I^I-b_I^Eb_E^I}.
\nonumber
\end{equation}
This limit is calculated using the limit 
of $N_I'(N_E)$ and proceeding in a similar way than in \cite{CCP}.
\item  If $b_E^Ib_I^E>b_E^Eb_I^I$ and 
$\frac{b_E^I}{V_F-V_R+b_I^I}
>\frac{b_E^E}{b_I^E}$, then 
$\displaystyle \lim_{N_E\to \infty}\mathcal{F}(N_E)=\infty$.

In this case the limit of $\mathcal{F}$ is a product of increasing
 functions with limit $\infty$.
\end{enumerate}

 \end{enumerate}

\emph{Step 4. The proof of the theorem is a consequence of the previous
steps.} 

The  monotonicity of  $\mathcal{F}$ and its limit, calculated in step 3, 
provide the number of steady states in terms of \eqref{condicion21}
and \eqref{condicion21bis},  since these conditions give the
range of the parameter values for which the limit of $\mathcal{F}$
can be compared to 1.

Under assumptions \eqref{condicion21} and \eqref{condicion3}, there
are no steady states. The reason is that  
for these values of $b_E^E$ the function $I_1(N_E)$ decreases, as
$b_E^Ib_I^E<b_E^Eb_I^I$ (see property 4.1 of $I_1(N_E)$). Thus,
\begin{equation}
I_1(N_E)<I_1(0)<\frac{b_E^E}{2(V_F+b_I^E N_I(\bar N_E))}
=\frac{1}{\bar N_E}<\frac{1}{N_E} \quad \forall \ N_E < \bar N_E.
\nonumber
\end{equation}
On the other hand, in an analogous way as in   
\cite{CCP}[Theorem 3.1(iv)], we can prove
\begin{equation}
I_1(N_E)< \frac{V_F-V_R}{b_E^EN_E-V_F-b_I^EN_I(N_E)}<\frac{1}{N_E}
 \qquad \forall \ N_E\ge\bar N_E,
\nonumber
\end{equation}
if  
\begin{equation}
w_F<0 \quad \mbox{and} \quad
\frac{V_F-V_R}{b_E^EN_E-V_F-b_I^EN_I(N_E)}<\frac{1}{N_E}
\quad \mbox{for all}   \quad N_E \ge \bar N_E.
\label{-}
\end{equation}
Therefore, to conclude this part of the proof 
we show \eqref{-}. Defining $g(N_E):= \frac{2[V_F+b_I^EN_I(N_E)]}{N_E}$,
since $b_E^E>\displaystyle \lim_{N_E \to \infty} g(N_E)=
\frac{2b_I^Eb_E^I}{V_F-V_R+b_I^I}$ and $g(0)=\infty$, 
there is  $N_E^*\ge0$ so that $g(N_E^*)=b_E^E$. 
As the monotonicity of $g$ is not known,  there might 
be several values of $N_E^*$ that solve the 
equation $g(N_E^*)=b_E^E$. However, the largest value $\bar{N}_E$ 
(see \eqref{eleccion_N}) ensures that $g(N_E)<b_E^E$ for all
 $N_E>\bar{N}_E$. Thus, 
$
b_E^E>g(N_E)> \frac{V_F+b_I^EN_I(N_E)}{N_E}$  for $N_E\ge \bar N_E$
and we obtain that  $w_F<0$  for $N_E \ge \bar N_E$.
The inequality $\frac{V_F-V_R}{b_E^EN_E-V_F-b_I^EN_I(N_E)}<\frac{1}{N_E}$,  
for all  $N_E \ge \bar N_E$, is proved using that 
 $b_E^E>2(V_F-V_R)$ and  $N_E>\frac{2[V_F+b_I^E(N_I(N_E)]}{b_E^E}$ due to
\eqref{condicion3} and \eqref{eleccion_N}.

To conclude the proof
we note that there is a unique steady state for parameters  where 
$\displaystyle \lim_{N_E \to \infty} \mathcal{F}(N_E)=\infty$.
Indeed, if $b_E^E$ is small enough, $\mathcal{F}$ is an increasing function. 
And, for $b_E^E$ small but such that the limit of
 $\mathcal{F}$ is finite, we deduce that
there is a unique stationary solution in an analogous way as 
 in \cite{CCP}[Theorem 3.1(i)] for a purely excitatory network.
\qed
As  proved in \cite{CCP} we can find conditions for the
connectivity parameters in order to have at least two steady states. 
We explain it in the following remark.
\begin{remark}
 If the parameters of the model satisfy \eqref{condicion21}
and 
\begin{eqnarray}
2a_Eb_E^E<\left[V_R+b_I^EN_I\left(\frac{2a_E}{(V_F-V_R)^2}\right)
\right](V_F-V_R)^2,
\label{condicion22}
\end{eqnarray}
there are at least two stationary solutions to \eqref{modelo}-\eqref{N}.

The proof is as follows.
As $\mathcal{F}(0)=0$ and   $\displaystyle 
\lim_{N_E \to \infty} \mathcal{F}(N_E)<1$
(due to \eqref{condicion21}), we have to prove that 
$\mathcal{F}(N_E)>1$  for $N_E \in D:=\left(
\frac{2a_E}{(V_F-V_R)^2},
\frac{V_R+b_I^EN_I\left(\frac{2a_E}{(V_F-V_R)^2}\right)}{b_E^E}\right)$.
This interval is not empty since \eqref{condicion22} holds,
and for all $N_E\in D$  we have
$
N_E \leq \frac{V_R+b_I^EN_I(N_E)}{b_E^E}$
because $N_I(N_E)$ is an increasing function.
Thus, $w_R>0$ for $N_E\in D$. 
Finally, in an analogous way as in 
case (ii) of Theorem 3.1 in \cite{CCP}, using \eqref{ecN2}
we obtain  
\begin{eqnarray*}
I_1(N_E) &\ge& \int_{w_R}^{w_F}\left[e^{-\frac{z^2}{2}}\int_{\textrm{max}(z,w_R)}^{w_F}e^{\frac {u^2}{2}}\ du\right]\ dz \\
&
=
&  
\int_{w_R}^{w_F}\left[e^{-\frac{z^2}{2}}
\int_{z}^{w_F}e^{\frac {u^2}{2}}\ du\right]\ dz\ge
\int_{w_R}^{w_F} \left[ \int_{
z
}^{w_F}  \ du \right]\ dz 
=
\frac{(V_F-V_R)^2}{2a_E}>\frac{1}{N_E}
\end{eqnarray*}
for  $N_E \in D$. 
\qed
\end{remark}

From a biological point of view, the previous analysis of the number of 
steady states shows the complexity of the set of stationary solutions
in terms of the model parameters: 
the reset and threshold potentials, $V_R$ and $V_F$, and the connectivity 
parameters, $b_E^E$, $b_E^I$, $b_I^E$ and $b_I^I$,
which describe 
the strengths of the synapses between excitatory and inhibitory neurons.
For example, in function of the connectivity parameter $b_E^E$, 
considering the rest of parameters fixed, we observe that if it is
small (there are  weak connections between excitatory neurons) there
is a unique steady state. Whereas, if it is large (strong 
connections between excitatory neurons) there are no steady states.
For intermediate values of $b_E^E$ different situations can occur: 
one, two or  three steady states.
In Section \ref{sec:numericalresults} we illustrate this complexity
with numerical experiments (see Figs. \ref{caso1}, \ref{caso2},
\ref{cruzados_chicos}, \ref{bifurcacion}).

\subsection{Long time behaviour}
In \cite{CCP} it was proved exponential convergence to the steady state for 
the linear case.
Later,  these results 
were extended in \cite{carrillo2014qualitative} to the nonlinear case,
for a purely excitatory
or inhibitory network
 with small connectivity parameters.
In both papers the main tools used were the relative entropy and Poincar\'e
inequalities. These techniques can be adapted for a  coupled
excitatory-inhibitory network, where the diffusion
terms are considered constant, $a_\alpha$, where $\alpha=E,I$. 
This is the goal of this 
subsection.

In Theorem \ref{th_eq} it  was proved that for small connectivity parameters
and constant diffusion terms
there is only one stationary solution of the system \eqref{modelo}-\eqref{N},
$\left(\rho_E^\infty, \rho_I^\infty, N_E^\infty,  N_I^\infty\right)$.
Therefore,  for any smooth convex function $G:\R^+\to \R$, we can define
the {\em relative entropy} for $\alpha= E, I$ as
\begin{equation}
\int_{-\infty}^{V_F} 
\rho_\alpha^\infty(v) G\left(\frac{\rho_\alpha(v,t)}{\rho_\alpha^\infty(v)} \right)
\, dv.
\nonumber
\end{equation}
Following similar computations to those developed in \cite{CCP, CP, 
carrillo2014qualitative} for fast-decaying solutions
to system \eqref{modelo}-\eqref{N}, i.e. 
weak solutions for which the weak formulation holds for all
test functions growing algebraically in $v$,
 it can be obtained,  for $\alpha=E,I$:
\begin{eqnarray}
\frac{d}{dt}
\int_{-\infty}^{V_F} 
\rho_\alpha^\infty(v) G\left(\frac{\rho_\alpha(v,t)}{\rho_\alpha^\infty(v)} \right)
\, dv  =  -  a_\alpha \int_{-\infty}^{V_F} \rho_\alpha^\infty(v) 
G''\left( \frac{\rho_\alpha(v,t)}{\rho_\alpha^\infty(v)} \right)
\left(\frac{\partial}{\partial v}
\frac{\rho_\alpha(v,t)}{\rho_\alpha^\infty(v)}\right)^2
\ dv
\nonumber
\\
 -  N_\alpha^\infty \left[
G\left(\frac{N_\alpha(t)}{N_\alpha^\infty}\right)-
G\left(\frac{\rho_\alpha(V_R,t)}{\rho_\alpha^\infty(V_R)}\right)
-\left(\frac{N_\alpha(t)}{N_\alpha^\infty}- 
\frac{\rho_\alpha(V_R,t)}{\rho_\alpha^\infty(V_R)}
 \right) G'\left(\frac{\rho_\alpha(V_R,t)}{\rho_\alpha^\infty(V_R)}\right)
\right]
\nonumber
\\
 +  \left(b_E^\alpha \bar{N}_E(t)-b_I^\alpha \bar{N}_I(t) \right)
\int_{-\infty}^{V_F} \frac{\partial\rho_\alpha^\infty(t)}{\partial v}
\left[G\left(\frac{\rho_\alpha(v,t)}{\rho_\alpha^\infty(v)}\right)
-\left(\frac{\rho_\alpha(v,t)}{\rho_\alpha^\infty(v)}\right)
G'\left(\frac{\rho_\alpha(v,t)}{\rho_\alpha^\infty(v)}\right)
   \right],
\label{entropiaalfa}
\end{eqnarray}
where $\bar{N}_\alpha(t)=N_\alpha(t)-N_\alpha^\infty$.
Therefore,  the  entropy production of the relative entropy
(its time derivative) for the excitatory  (resp. inhibitory) population
depends on the firing rate of the inhibitory (resp. excitatory) population.
In other words, both entropy productions are linked by means of the firing
rates. However, for the quadratic entropy, $G(x)=\left(x-1\right)^2$,
 a control on the sum of them,
\begin{equation}
E[t]:=\int_{-\infty}^{V_F} 
\left[
\rho_E^\infty(v) 
\left(\frac{\rho_E(v,t)}{\rho_E^\infty(v)} -1 \right)^2
+
\rho_I^\infty(v) \left(\frac{\rho_I(v,t)}{\rho_I^\infty(v)} -1 \right)^2
\right]  \, dv,
\label{entropy}
\end{equation}
 can be obtained for small connectivity parameters,  
  proving an analogous result as 
in \cite{carrillo2014qualitative}[Theorem 2.1].  
\begin{theorem}
Assume $a_\alpha$ constant for $\alpha=E,I$, the connectivity parameters $b_E^E$, $b_E^I$, $b_I^E$ and $b_I^I$ 
small enough and  initial data $(\rho_E^0, \rho_I^0)$ such that
\begin{equation}
E[0]
<\frac{1}{2 \max\left(b_E^E+b_I^E,b_E^I+b_I^I \right)}.
\label{hip_medio}
\end{equation}
Then, for fast decaying solutions to \eqref{modelo}-\eqref{N} 
there is a constant $\mu>0$ such that, for all $t\ge 0$, 
$$
E[t] \le e^{-\mu t} \, E[0].
$$
Consequently, for $\alpha=E,I$
$$
\int_{-\infty}^{V_F} \rho_\alpha^\infty \left(
\frac{\rho_\alpha(v,t)}{\rho_\alpha^\infty(v)} -1\right)^2\, dv
\le e^{-\mu t} \, E[0].
$$
\end{theorem}  
\proof 
The proof follows analogous steps as  in 
\cite{carrillo2014qualitative}[Theorem 2.1], with the main 
difficulty that, in this case, the entropy productions
for excitatory and inhibitory populations are linked. 
This is the reason why  the total relative entropy, given by the 
functional \eqref{entropy}, has to be considered. 
Thus, the entropy production is the
 sum of the entropy productions of each population.
In this way, the terms of $\bar{N}_\alpha$ can be gathered and bound properly,
as it is shown below.

Using \eqref{entropiaalfa} with $G(x)=(x-1)^2$   
for each term of the entropy \eqref{entropy},  
its time derivative can be written as 
\begin{eqnarray}
\frac{d}{dt} E[t] & =& 
  2a_E \int_{-\infty}^{V_F} \rho_E^\infty(v) 
\left(\frac{\partial}{\partial v}
\frac{\rho_E(v,t)}{\rho_E^\infty(v)}\right)^2 \ dv
 -  2a_I \int_{-\infty}^{V_F} \rho_I^\infty(v) 
\left(\frac{\partial}{\partial v}\frac{\rho_I(v,t)}{\rho_I^\infty(v)}\right)^2
\ dv
\nonumber
\\
& - &  N_E^\infty \left(\frac{N_E(t)}{N_E^\infty}-
\frac{\rho_E(V_R,t)}{\rho_E^\infty(V_R)}\right)^2
 -  N_I^\infty \left(\frac{N_I(t)}{N_I^\infty}-
\frac{\rho_I(V_R,t)}{\rho_I^\infty(V_R)}\right)^2
\label{derivada_entropia}
\\
& + & 2\left(b_E^E \bar{N}_E(t)-b_I^E \bar{N}_I(t) \right)
\int_{-\infty}^{V_F} \rho_E^\infty(v)\left[\frac{\partial}{\partial v}
\frac{\rho_E(v,t)}{\rho_E^\infty(v)}
\left(\frac{\rho_E(v,t)}{\rho_E^\infty(v)} -1\right)+
\frac{\partial}{\partial v}\frac{\rho_E(v,t)}{\rho_E^\infty(v)}\right] \ dv
\nonumber
\\
& + &  2\left(b_E^I \bar{N}_E(t)-b_I^I \bar{N}_I(t) \right)
\int_{-\infty}^{V_F} \rho_I^\infty(v)\left[\frac{\partial}{\partial v}
\frac{\rho_I(v,t)}{\rho_I^\infty(v)}
\left(\frac{\rho_I(v,t)}{\rho_I^\infty(v)} -1\right)+
\frac{\partial}{\partial v}\frac{\rho_I(v,t)}{\rho_I^\infty(v)}\right] \ dv,
\nonumber
\end{eqnarray}
where the last two terms were obtained using that $G(x)-xG'(x)=1-x^2$, 
integrating by parts and taking into account the boundary conditions 
\eqref{ec4}.

Considering the inequality $(a+b)^2\ge \epsilon \left(a^2-2b^2\right)$, for
$a, \, b \in \R$ and $0<\epsilon <1/2$, the Sobolev injection of 
$L^\infty(I)$ in 
$H^1(I)$ for a small neighborhood $I$ of 
$V_R$, where $\rho_\alpha^\infty$ is bounded from below
and the Poincar\'e inequality (see \cite{CCP}[Appendix] and
\cite{carrillo2014qualitative} for details), the third and fourth terms
in  \eqref{derivada_entropia} satisfy:
\begin{eqnarray*}
 -  N_\alpha^\infty \left(\frac{N_\alpha(t)}{N_\alpha^\infty}-
\frac{\rho_\alpha(V_R,t)}{\rho_\alpha^\infty(V_R)}\right)^2 
\leq -C_0^\alpha  \left(\frac{N_\alpha(t)}{N_\alpha^\infty} -1 \right)^2
+\frac{a_\alpha}{2}\int_{-\infty}^{V_F}\rho_\alpha^\infty(v)
\left( \frac{\partial}{\partial v} 
\frac{\rho_\alpha(v,t)}{\rho_\alpha^\infty}(v,t) \right) ^2dv,
\end{eqnarray*}
for some positive constant $C_0^\alpha$.

The estimate of the last two terms of 
 \eqref{derivada_entropia} is quite more involved. 
First, each  term is split into four addends, so that
\begin{eqnarray}
2\left(b_E^\alpha \bar{N}_E(t)-b_I^\alpha \bar{N}_I(t) \right) 
\int_{-\infty}^{V_F} \rho_\alpha^\infty(v) & & 
\!  \!  \!  \!  \!  \!  \!  \!  \!  \!  \!  \!  \!  \!  
\left[\frac{\partial}{\partial v}
\frac{\rho_\alpha(v,t)}{\rho_\alpha^\infty(v)}
\left(\frac{\rho_\alpha(v,t)}{\rho_\alpha^\infty(v)} -1\right)+
\frac{\partial}{\partial v}\frac{\rho_\alpha(v,t)}{\rho_\alpha^\infty(v)}\right] \ dv
 \nonumber
 \\
& = & 2b_E^\alpha \bar{N}_E(t)\int_{-\infty}^{V_F} 
\rho_\alpha^\infty(v)\frac{\partial}{\partial v}\frac{\rho_\alpha(v,t)}{\rho_\alpha^\infty(v)}
\left(\frac{\rho_\alpha(v,t)}{\rho_\alpha^\infty(v)} -1\right) \ dv
 \nonumber
 \\
&  & +  2b_E^\alpha \bar{N}_E(t) \int_{-\infty}^{V_F} 
\rho_\alpha^\infty(v)\frac{\partial}{\partial v}
\frac{\rho_\alpha(v,t)}{\rho_\alpha^\infty(v)} \ dv
 \nonumber
 \\
& & - 2b_I^\alpha \bar{N}_I(t)\int_{-\infty}^{V_F} \rho_\alpha^\infty(v)
\frac{\partial}{\partial v}\frac{\rho_\alpha(v,t)}{\rho_\alpha^\infty(v)}
\left(\frac{\rho_\alpha(v,t)}{\rho_\alpha^\infty(v)} -1\right) \ dv
 \nonumber
\\
&  & -   2b_I^\alpha \bar{N}_I(t)\int_{-\infty}^{V_F} 
\rho_\alpha^\infty(v)\frac{\partial}{\partial v}
\frac{\rho_\alpha(v,t)}{\rho_\alpha^\infty(v)} \ dv. \nonumber
\end{eqnarray}
Then, Cauchy-Schwarz's and Young's inequalities provide
\begin{eqnarray*}
2b_E^\alpha & & \!  \!  \!  \!  \!  \!  \!  \!  \!  \!  \!  \!  \!  \!  
   |\bar{N}_E(t)| 
\int_{-\infty}^{V_F}\rho_\alpha^\infty(v)\left|\frac{\partial}{\partial v}
\frac{\rho_\alpha(v,t)}{\rho_\alpha^\infty(v)}
\left(\frac{\rho_\alpha(v,t)}{\rho_\alpha^\infty(v)} -1\right)\right| \, dv 
\\
& \leq & 
 \frac{b_E^\alpha}{a_E}{(N_E^\infty)}^2
\left(\frac{N_E(t)}{N_E^\infty}-1\right)^2
 + a_E b_E^\alpha \int_{-\infty}^{V_F} \rho_\alpha^\infty(v)
\left(\frac{\partial}{\partial v}
\frac{\rho_\alpha(v,t)}{\rho_\alpha^\infty(v)}\right)^2 \, dv
\int_{-\infty}^{V_F} \rho_\alpha^\infty(v)
\left(\frac{\rho_\alpha(v,t)}{\rho_\alpha^\infty(v)} -1\right)^2\, dv 
\end{eqnarray*}
and
\begin{eqnarray*}
2b_E^\alpha |\bar{N}_E(t)| 
\int_{-\infty}^{V_F} \rho_\alpha^\infty(v)
\left|\frac{\partial}{\partial v}
\frac{\rho_\alpha(v,t)}{\rho_\alpha^\infty(v)}\right| \, dv 
\leq  
4\frac{{(b_E^\alpha)}^2}{a_E}{N_E^\infty}^2
\left(\frac{N_E(t)}{N_E^\infty}-1 \right)^2+\frac{a_E}{4} 
\int_{-\infty}^{V_F} \rho_\alpha^\infty(v)\left(\frac{\partial}{\partial v}
\frac{\rho_\alpha(v,t)}{\rho_\alpha^\infty(v)}\right)^2 \, dv.
\end{eqnarray*}
Getting  these bounds together
\begin{eqnarray}
\frac{d}{dt}
E[t]
& \leq &  
\left( \frac{b_E^E}{a_E}+   
\frac{b_E^I}{a_I} + 4 \frac{{(b_E^E)}^2}{a_E}+ 
4\frac{{(b_E^I)}^2}{a_I}-\frac{C_0^E}{{N_E^\infty}^2} \right){(N_E^\infty)}^2
\left(\frac{N_E(t)}{N_E^\infty}-1\right)^2  \label{cota}
\\
& &+   \left( \frac{b_I^E}{a_E}+   \frac{b_I^I}{a_I} 
+ 
4 \frac{{(b_I^E)}^2}{a_E}+4 \frac{{(b_I^I)}^2}{a_I} 
-\frac{C_0^I}{{N_I^\infty}^2}\right){(N_I^\infty)}^2
\left(\frac{N_I(t)}{N_I^\infty}-1\right)^2 
\nonumber
\\
& & 
 -  a_E\int_{-\infty}^{V_F} \rho_E^\infty(v)
\left( \frac{\partial}{\partial v} \frac{\rho_E(v,t)}{\rho_E^\infty}(v,t) 
\right)^2 dv \left[ 1-(b_E^E+b_E^I)
\int_{-\infty}^{V_F} \rho_E^\infty(v)\left(  
\frac{\rho_E(v,t)}{\rho_E^\infty}(v,t)-1 \right)^2dv \right] 
\nonumber
\\
& & 
- a_I\int_{-\infty}^{V_F} \rho_I^\infty(v)
\left( \frac{\partial}{\partial v} 
\frac{\rho_I(v,t)}{\rho_I^\infty}(v,t) \right)^2dv 
\left[ 1-(b_I^E+b_I^I)
\int_{-\infty}^{V_F} \rho_I^\infty(v)\left(
\frac{\rho_I(v,t)}{\rho_I^\infty}(v,t)-1 \right)^2dv \right].
\nonumber
\end{eqnarray}
In this way, for  $b_k^\alpha$
small enough such that
$
\left( \frac{b_\alpha^E}{a_E}+   
\frac{b_\alpha^I}{a_I} + 
4 \frac{{(b_\alpha^E)}^2}{a_E}+ 
\frac{{(b_\alpha^I)}^2}{a_I}-
\frac{C_0^\alpha}{{(N_\alpha^\infty)}^2} \right)<0,
$
the first and second terms of the right hand side 
of \eqref{cota} are negative, thus 
\begin{eqnarray}
\frac{d}{dt}  E[t] 
&\!\! \leq &  \! \!
-  a_E\int_{-\infty}^{V_F} \rho_E^\infty(v)\left( 
\frac{\partial}{\partial v} 
\frac{\rho_E(v,t)}{\rho_E^\infty}(v,t) \right) ^2 \!\!\! \!dv 
\left[ 1-(b_E^E+b_E^I)\int_{-\infty}^{V_F} 
\rho_E^\infty(v)\left(  \frac{\rho_E(v,t)}{\rho_E^\infty}(v,t)-1 \right) ^2
\!\!\! \!
dv \right] 
\nonumber
\\
& \! \!\!\!\! \! - & \! \!\!\! \! \!
a_I\int_{-\infty}^{V_F} \rho_I^\infty(v)
\left( \frac{\partial}{\partial v} 
\frac{\rho_I(v,t)}{\rho_I^\infty}(v,t) \right) ^2\!\!\! \!  dv 
\left[ 1-(b_I^E+b_I^I)\int_{-\infty}^{V_F} \rho_I^\infty(v)
\left(  \frac{\rho_I(v,t)}{\rho_I^\infty}(v,t)-1 \right) ^2
 \!\!\! \!  dv \right]\!\!.
\label{cota_mejorada}
\end{eqnarray}
Denoting $C=\max\left(b_E^E+b_I^E,b_E^I+b_I^I \right)$,
the entropy production can be bounded as follows
\begin{equation}
 \frac{d}{dt} E[t] \leq 
 - 
\int_{-\infty}^{V_F}
\left[
a_E \rho_E^\infty(v)\left( \frac{\partial}{\partial v} 
\frac{\rho_E(v,t)}{\rho_E^\infty}(v,t) \right) ^2
+a_I\rho_I^\infty(v)\left( \frac{\partial}{\partial v} 
\frac{\rho_I(v,t)}{\rho_I^\infty}(v,t) \right) ^2\right] dv 
\left( 1- C E[t] \right).
\label{cota_mejorada2}
\end{equation}
Now, using \eqref{hip_medio} and Gronwall's inequality 
it can be proved that $E[t]$ 
decreases for all times
and  $E[t]<\frac{1}{2 C}$,  for all $t>0$.
Thus,
$$
\frac{d}{dt} E[t]
\leq -\frac{a_E}{2}\int_{-\infty}^{V_F} 
\rho_E^\infty(v)\left( \frac{\partial}{\partial v} 
\frac{\rho_E(v,t)}{\rho_E^\infty}(v,t) \right) ^2dv-
\frac{a_I}{2}\int_{-\infty}^{V_F} 
\rho_I^\infty(v)\left( \frac{\partial}{\partial v} 
\frac{\rho_I(v,t)}{\rho_I^\infty}(v,t) \right) ^2dv.
$$
Applying the Poincar\'e inequality on each term, there exist
$\gamma, \, \gamma'>0$ such that
$$
\frac{d}{dt} E[t]
\leq -\frac{a_E \gamma}{2}\int_{-\infty}^{V_F} \rho_E^\infty(v) \left(
\frac{\rho_\alpha(v,t)}{\rho_E^\infty(v)} -1\right)^2\,dv-
\frac{a_I \gamma '}{2}\int_{-\infty}^{V_F} \rho_I^\infty(v) \left(
\frac{\rho_\alpha(v,t)}{\rho_I^\infty(v)} -1\right)^2\,dv.
$$
Considering now
 $\mu=\min\left(\frac{a_E \gamma}{2},\frac{a_I \gamma ' }{2}\right)$ we obtain
$\frac{d}{dt} E[t] \leq -\mu E[t]$.
Finally, Gronwall's inequality concludes the proof.
\qed

\section{Numerical results}
\label{sec:numericalresults}

The analytical results proved in previous sections are shown numerically  in the present section. The numerical scheme considered for this purpose is based on  a fifth order conservative finite difference WENO 
(Weighted Essentially Non-Oscillatory)
scheme for the advection term, standard second order centered finite differences for the diffusion  term and an explicit third order TVD  
(Total Variation Diminishing)
Runge-Kutta scheme for the time evolution. To reduce the computation time, 
parallel computation techniques for a two cores code is developed. 
Thus, the time evolution for both equations of the system is 
calculated simultaneously.
Each core handles one of the equations. 
MPI  
(Message Passing Interface)
communication between the cores has been included in the code, since the system is coupled by the firing rates. Therefore, at the end of each Runge-Kutta step, each core needs to know the value of the firing rate of the other core. 

For the simulations, an uniform mesh for  $v \in [-V_{left}, V_F]$ 
is considered. The value of $-V_{left}$ is chosen such that 
$\rho_\alpha(-V_{left},t) \sim 0$
(since $\rho_\alpha(-\infty,t)=0$).  
The time step size is adapted dynamically during the simulations via a 
CFL  
(Courant-Friedrich-Levy)
 time step condition.  Some parameter values are common to most simulations, $V_F=2$, $V_R=1$ and $ \nu_{E,ext}=0$ and the diffusion 
terms $a_\alpha(N_E,N_I)$  have been taken constant as $a_\alpha=1$. 
In the simulations where these values are different, the 
considered values are indicated in their figures and explanations.
In most cases, the initial condition is
\begin{equation}
\rho^0_\alpha(v)=
\frac{k}{\sqrt{2 \pi}}e^{-\frac{(v-v_0^\alpha)^2}{2 {\sigma_0^\alpha}^2}},
\label{ci_maxwel}
\end{equation}
where $k$ is a constant such that 
$\displaystyle\int_{-V_{left}}^{V_F}\rho^0_\alpha(v) \ dv \approx 1$ numerically.
However,  in order to analyze the stability of steady states, 
stationary profiles are taken as initial conditions  
\begin{equation}
\rho_\alpha^0(v) = 
\frac{N_\alpha}{a_\alpha(N_E,N_I)}
e^{-\frac{(v-V_0^\alpha(N_E,N_I))^2}{2a_\alpha(N_E,N_I)}}
\int_{\max(v,V_R)}^{V_F}
e^{\frac{(w-V_0^\alpha(N_E,N_I))^2}{2a_\alpha(N_E,N_I)}} \ dw, 
\quad \alpha=E,I, \label{soleq}
\end{equation}
with $V_0^\alpha (N_E,N_I)=b_E^\alpha N_E - b_I^\alpha N_I + 
(b_E^\alpha-b_E^E)\nu_{E,ext}$  and where $N_\alpha$ is an approximated value of the
 stationary firing rate.

\

\begin{figure}[H]
\begin{center}
\begin{minipage}[c]{0.33\linewidth}
\begin{center}
\includegraphics[width=\textwidth]{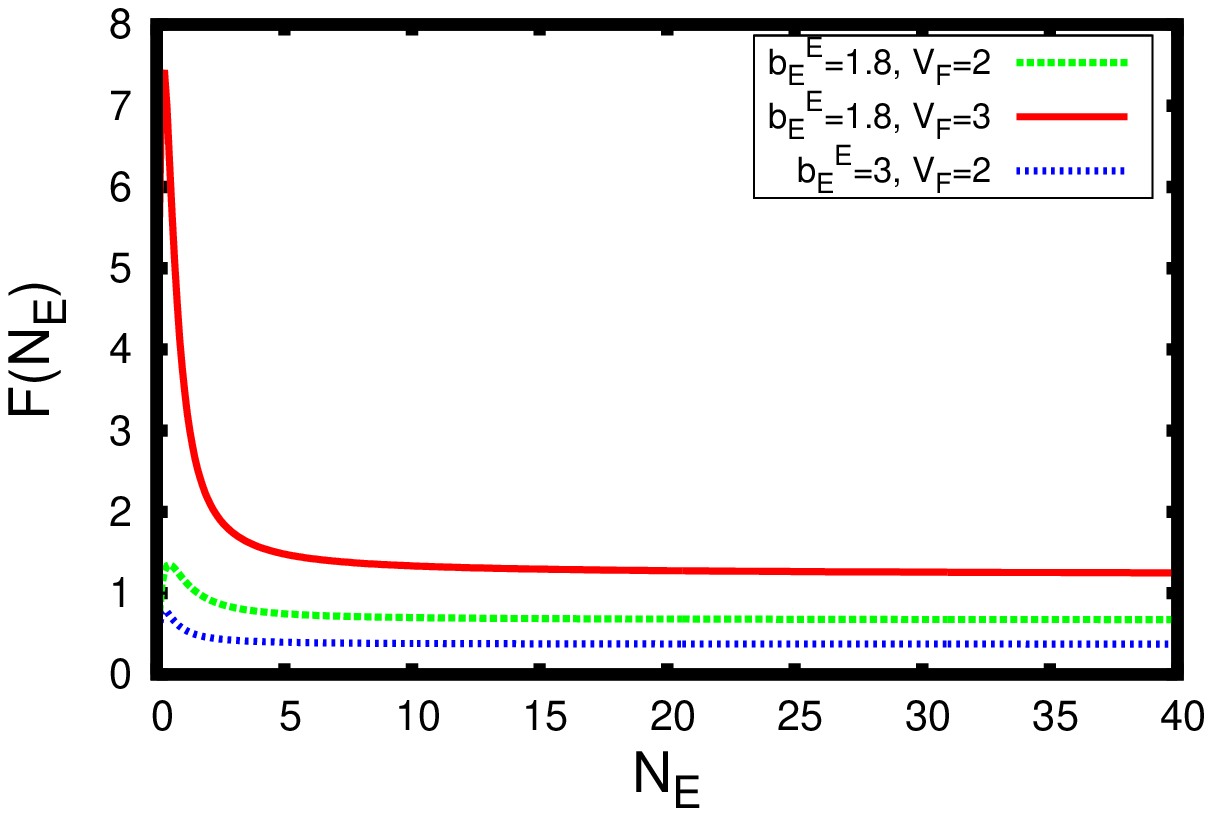}
\end{center}
\end{minipage}
\begin{minipage}[c]{0.33\linewidth}
\begin{center}
\includegraphics[width=\textwidth]{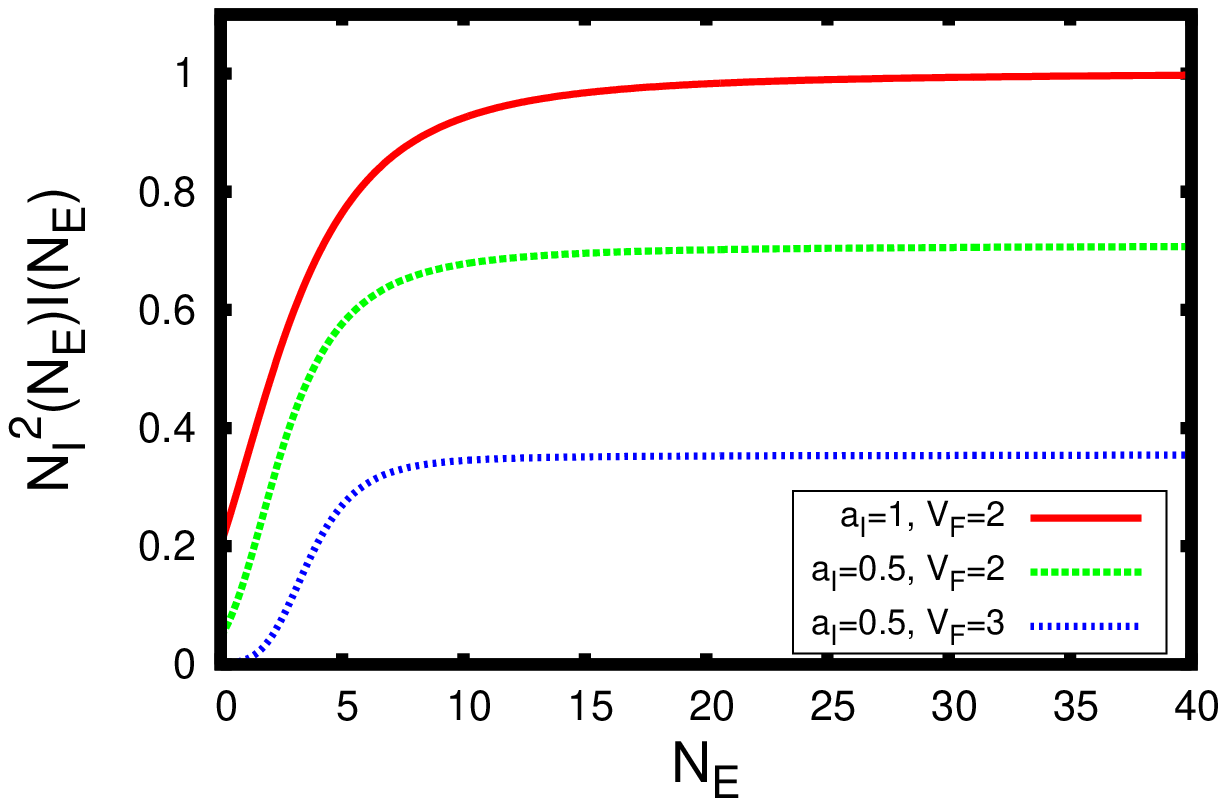}
\end{center}
\end{minipage}
\end{center}
\caption{
Study of the limits of $\mathcal{F}(N_E)$ and $N_I^2(N_E)I(N_E)$ (see
\eqref{N'bis}) when it is finite.
\newline
Left figure: 
$b_I^E=0.75$, $b_E^I=0.5$, $b_I^I=0.25$, $a_E=1$, $a_I=1$
 for different values of $b_E^E$ and $V_F$. 
\newline
Right figure: 
$b_E^E=1.8$, $b_I^E=0.75$, $b_E^I=0.5$, $b_I^I=0.25$ $a_E=1$  
for different values of $a_I$ and $V_F$.}
\label{limiteF}
\end{figure}

To get an idea about the number of steady states, 
the system \eqref{ecN2} is solved numerically. For 
every $N_E$ fixed, $N_I(N_E)$ is calculated as the root of 
$N_I(N_E)I_2(N_E,N_I(N_E))-1=0$ using the bisection method with 
tolerance $10^{-8}$, and then a numerical approximation 
of $\mathcal{F}(N_E)$ is developed 
by quadrature formulas
 to find the number
of intersections with the function 1. Fig. \ref{limiteF} 
shows the graphics of $N_I^2(N_E)I(N_E)$ (see
\eqref{N'bis}) and $\mathcal{F}(N_E)$ for parameter values for which
the limit of  $\mathcal{F}(N_E)$ is finite. The numerical approximation
of the limits is in agreement with their analytical limits.

\begin{figure}[H]
\begin{minipage}[c]{0.33\linewidth}
\begin{center}
\includegraphics[width=\textwidth]{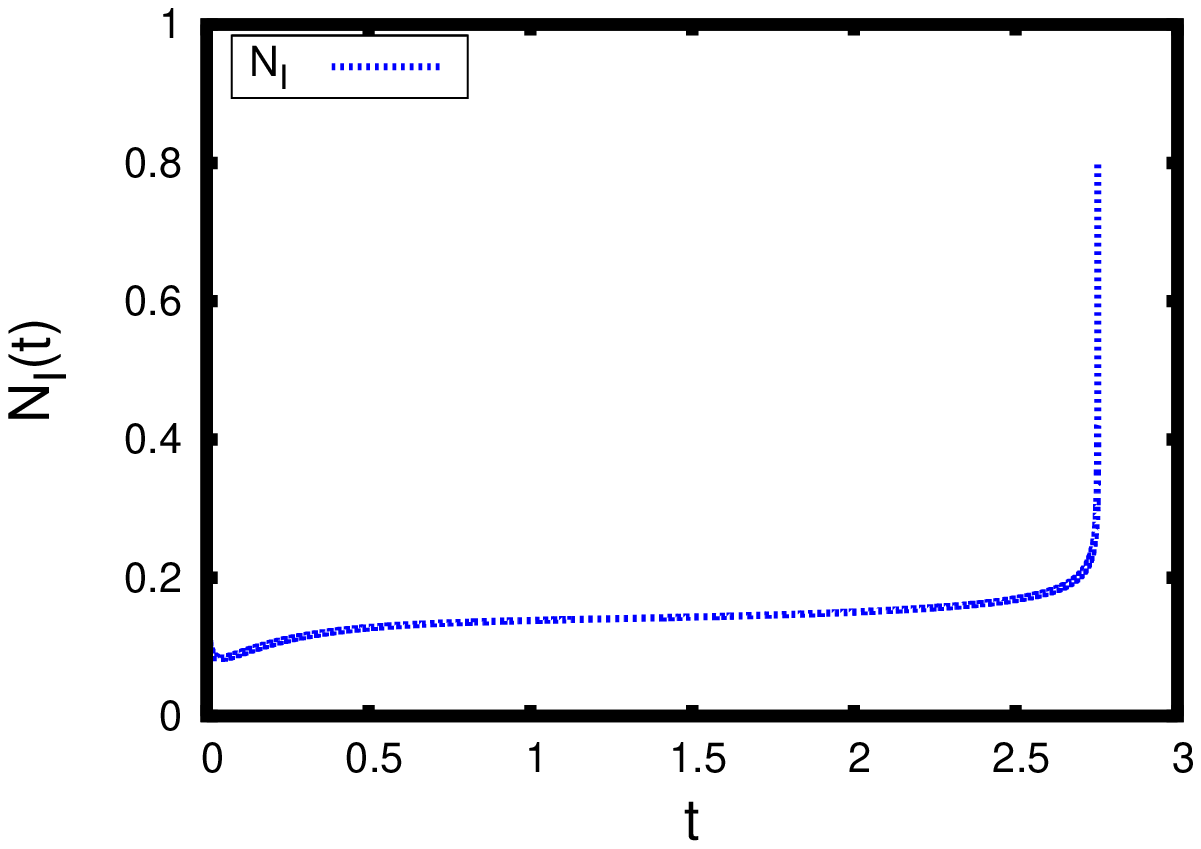}
\end{center}
\end{minipage}
\begin{minipage}[c]{0.33\linewidth}
\begin{center}
\includegraphics[width=\textwidth]{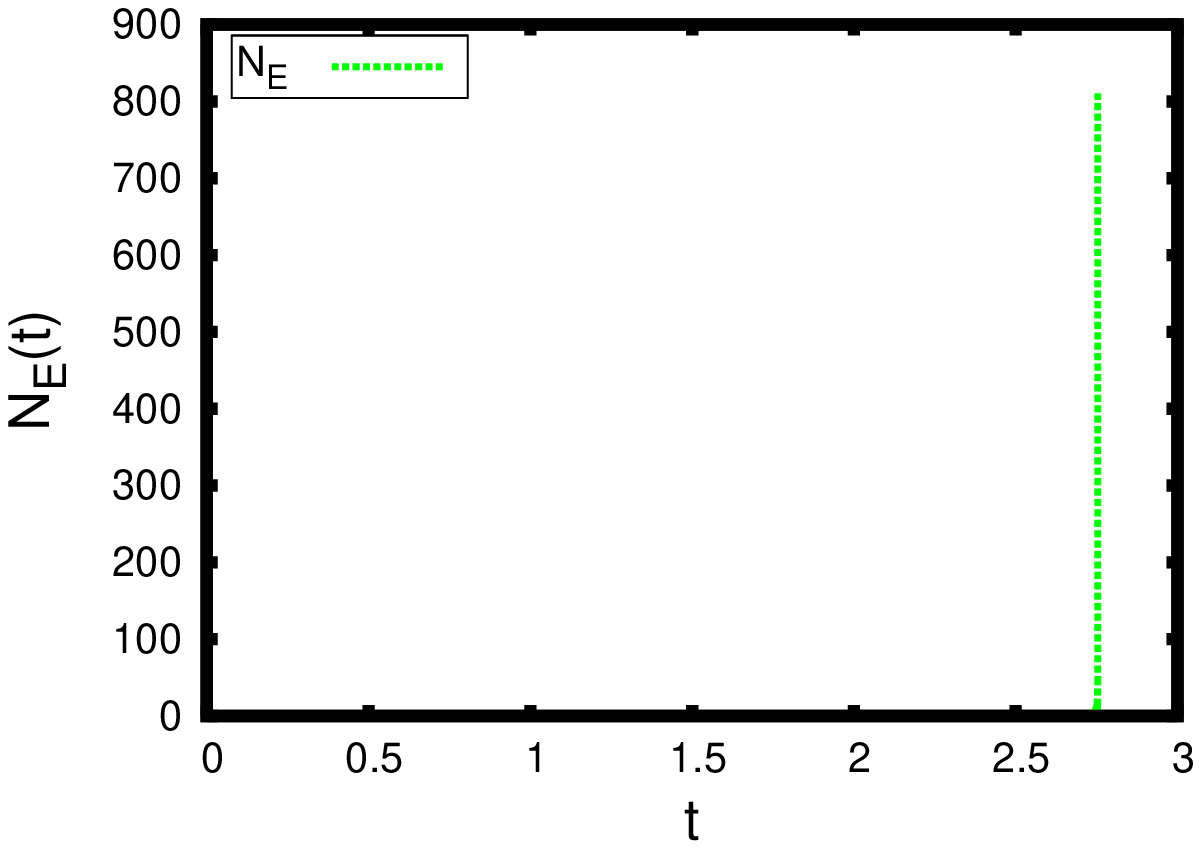}
\end{center}
\end{minipage}
\begin{minipage}[c]{0.33\linewidth}
\begin{center}
\includegraphics[width=\textwidth]{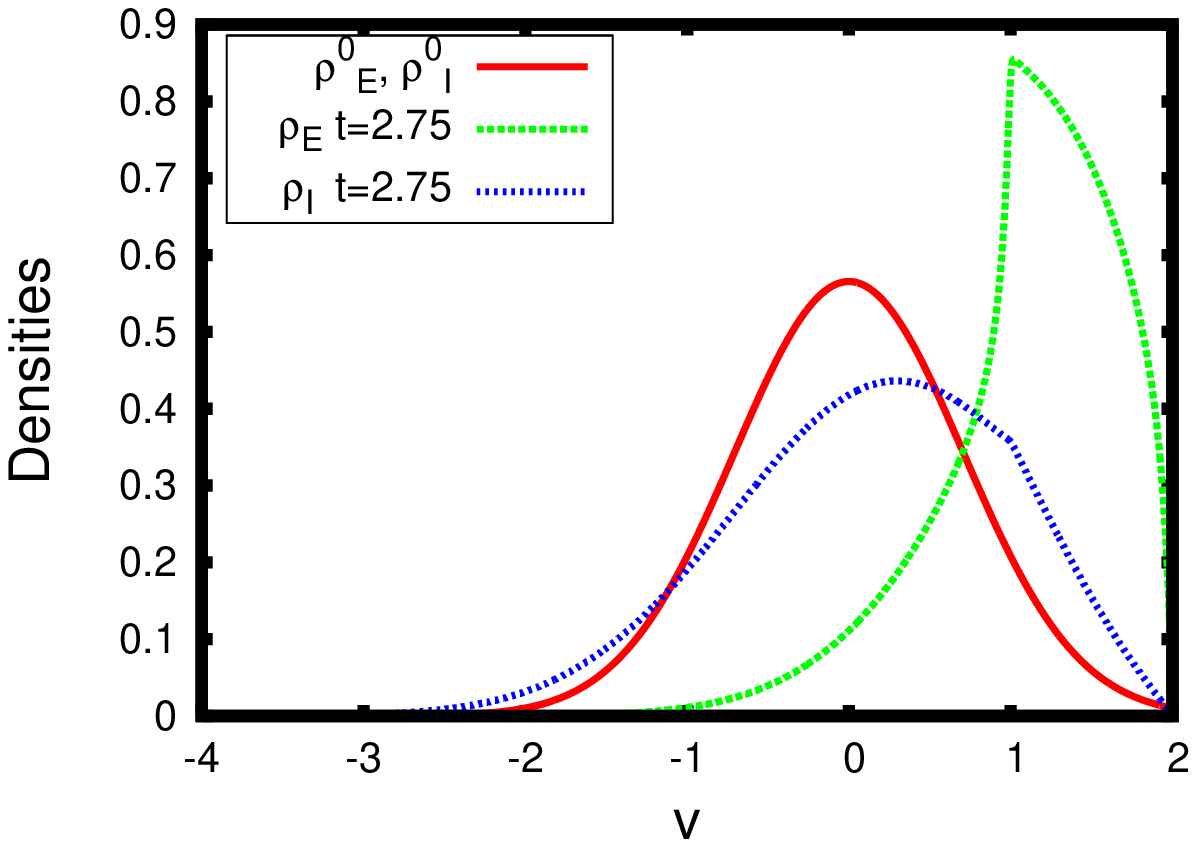}
\end{center}
\end{minipage}
\caption{Firing rates and probability densities for $b_E^E=3$, $b_I^E=0.75$, 
$b_E ^I=0.5$, $b_I^I=0.25$, in case of a normalized Maxwellian initial 
condition with mean  0 and variance 0.5 (see \eqref{ci_maxwel}). $N_E$ blows-up because of the large
value of $b_E ^E$.}
\label{blowup_bEE}
\end{figure}

\begin{figure}[H]
\begin{minipage}[c]{0.33\linewidth}
\begin{center}
\includegraphics[width=\textwidth]{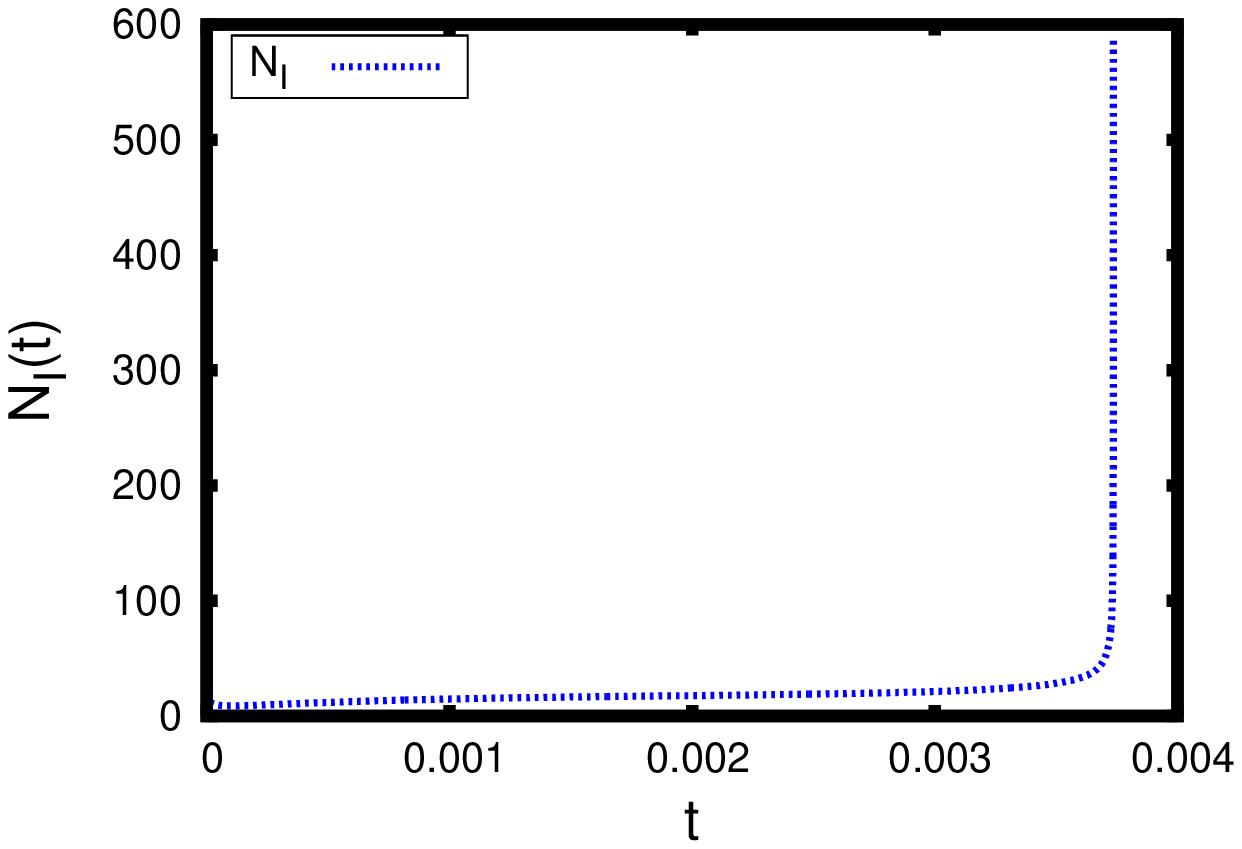}
\end{center}
\end{minipage}
\begin{minipage}[c]{0.33\linewidth}
\begin{center}
\includegraphics[width=\textwidth]{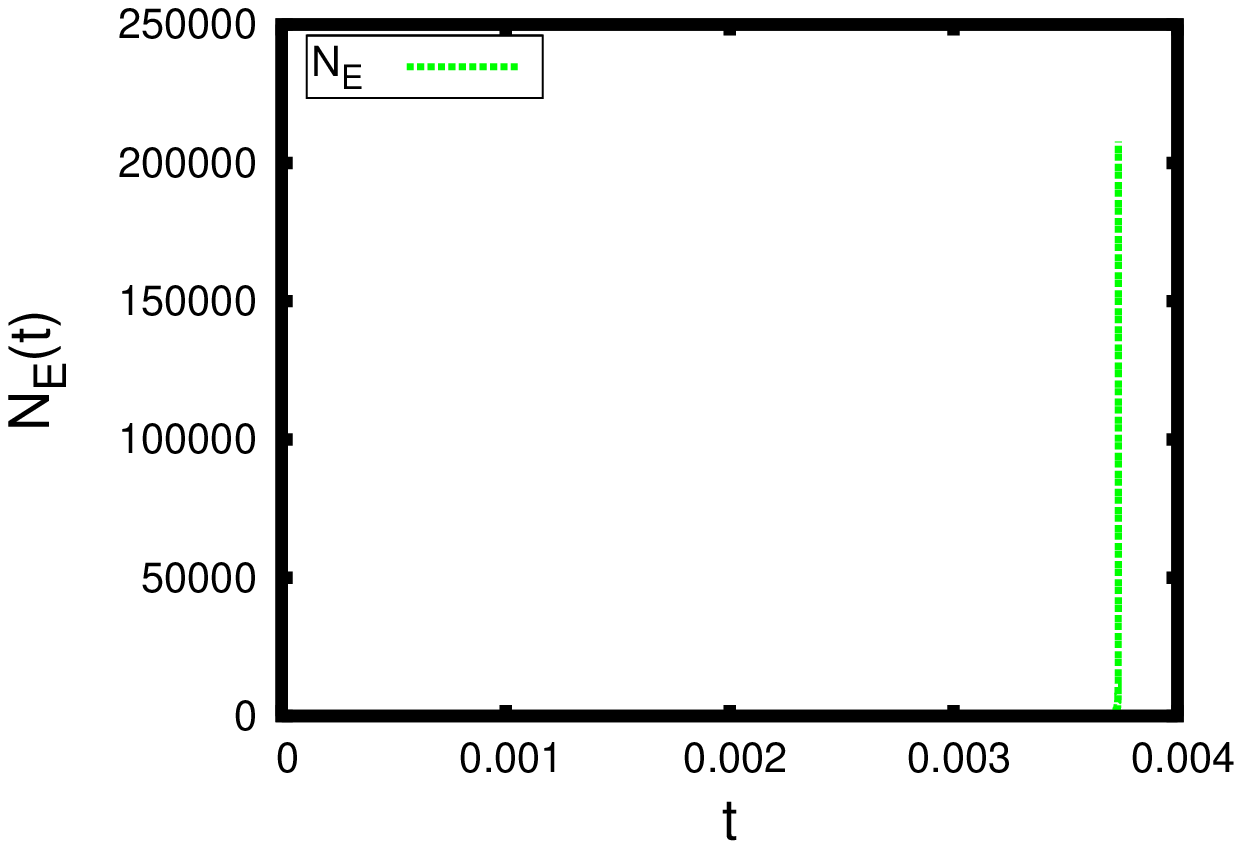}
\end{center}
\end{minipage}
\begin{minipage}[c]{0.33\linewidth}
\begin{center}
\includegraphics[width=\textwidth]{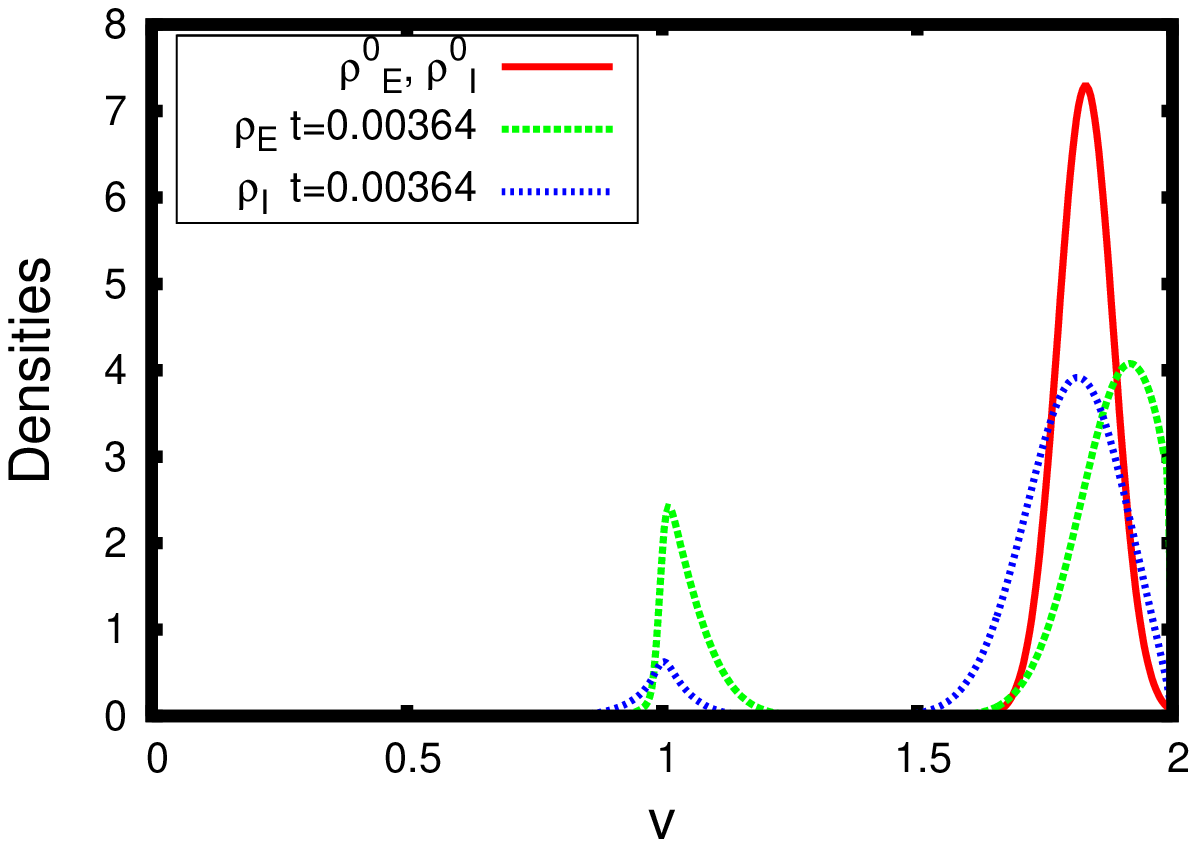}
\end{center}
\end{minipage}
\caption{Firing rates and probability densities for $b_E^E=0.5$, $b_I^E=0.25$, 
$b_E ^I=0.25$, $b_I^I=1$, in case of a normalized concentrated Maxwellian 
initial condition with mean  1.83 and variance 0.003 (see \eqref{ci_maxwel}).The initial condition concentrated
close to $V_F$ provokes the blow-up of $N_E$.}
\label{blowup_ci}
\end{figure}

\emph{\textbf{Blow-up.}} Numerical results for 
three situations in which the solutions  are not  global-in-time are described. 
In Fig. \ref{blowup_bEE} we depict  a blow-up situation produced by a 
big value of $b_E^E$, where the initial datum is far from $V_F$. 
However, for a value of the connectivity parameter $b_E^E$  small,
we show a blow-up situation originated  by an initial condition 
quite concentrated around $V_F$ in Fig. \ref{blowup_ci}.
In both cases we observe that  an extremely fast increase of the
firing rate of the excitatory population causes the blow-up
of the system. Furthermore, we notice that 
the firing rate of the inhibitory population also starts to grow sharply.
Nevertheless, the values it takes are quite small in comparison with 
those from the excitatory population.

\begin{figure}[H]
\begin{minipage}[c]{0.33\linewidth}
\begin{center}
\includegraphics[width=\textwidth]{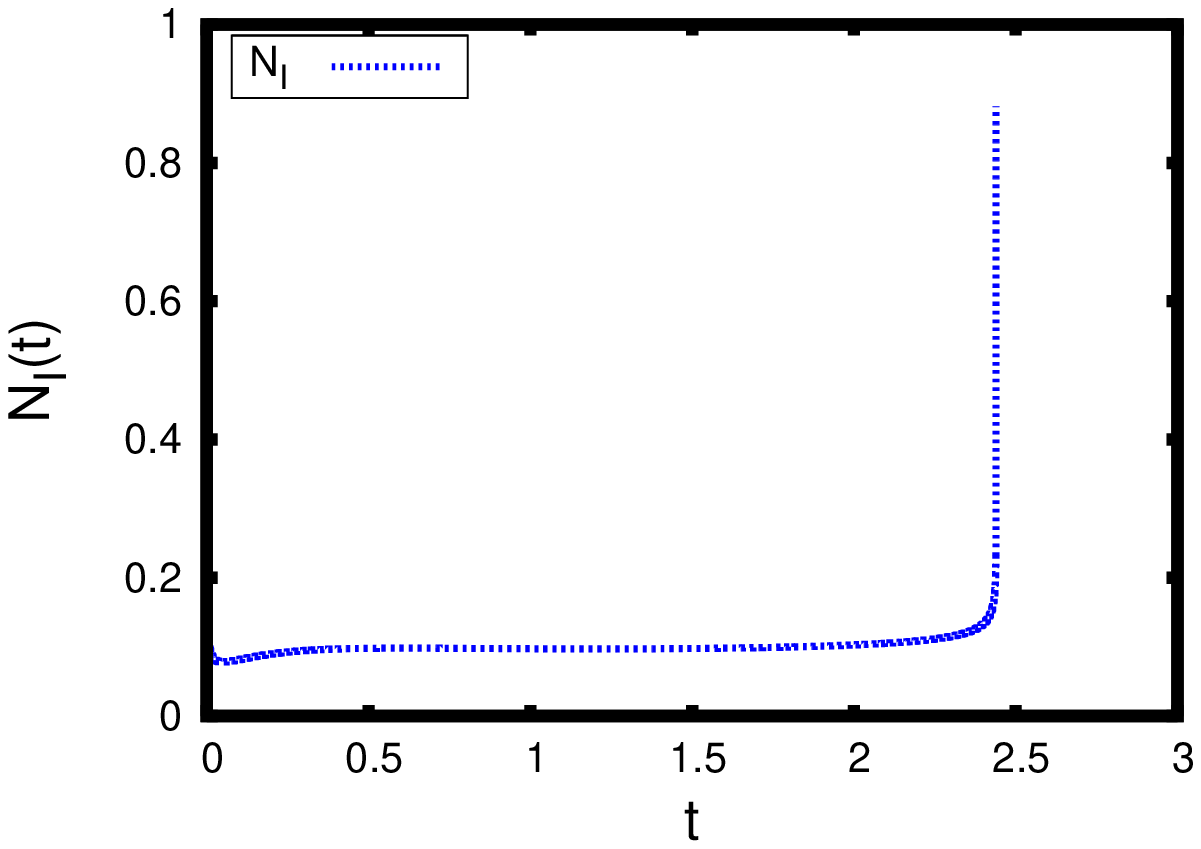}
\end{center}
\end{minipage}
\begin{minipage}[c]{0.33\linewidth}
\begin{center}
\includegraphics[width=\textwidth]{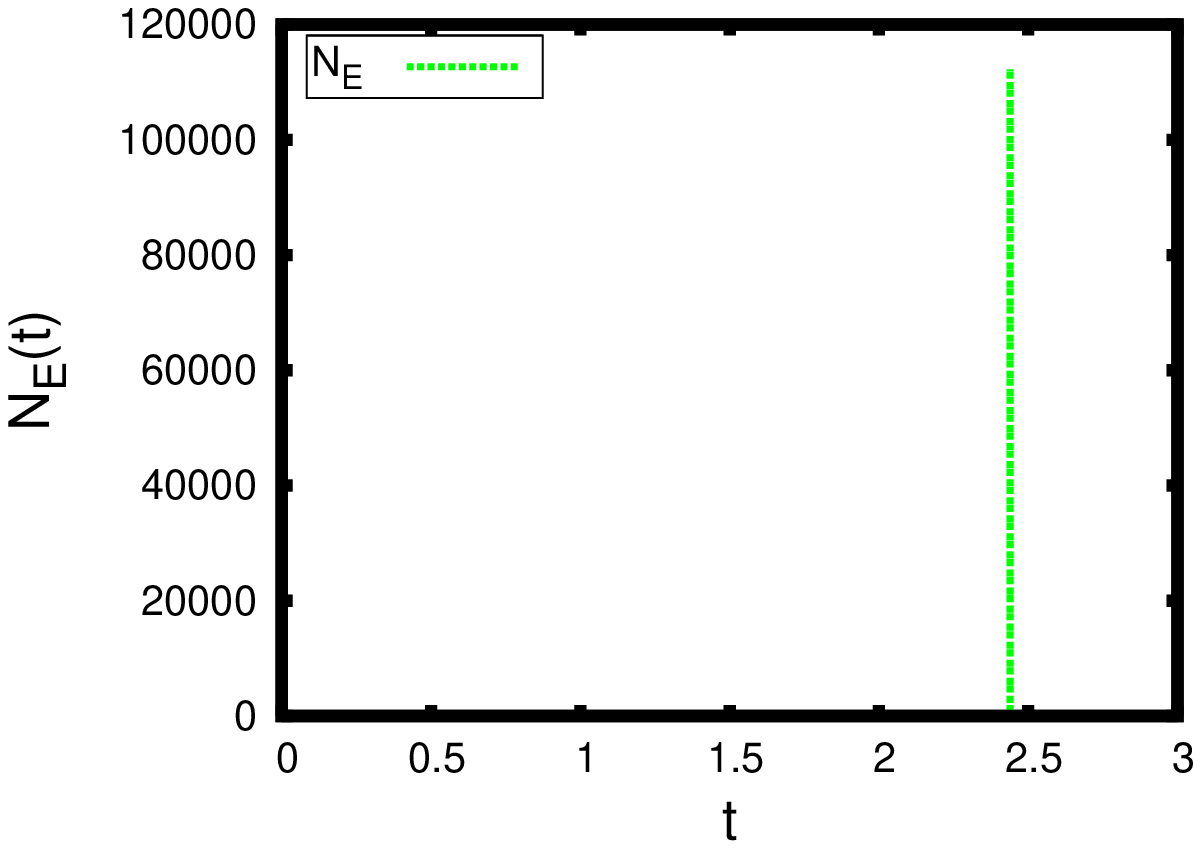}
\end{center}
\end{minipage}
\begin{minipage}[c]{0.33\linewidth}
\begin{center}
\includegraphics[width=\textwidth]{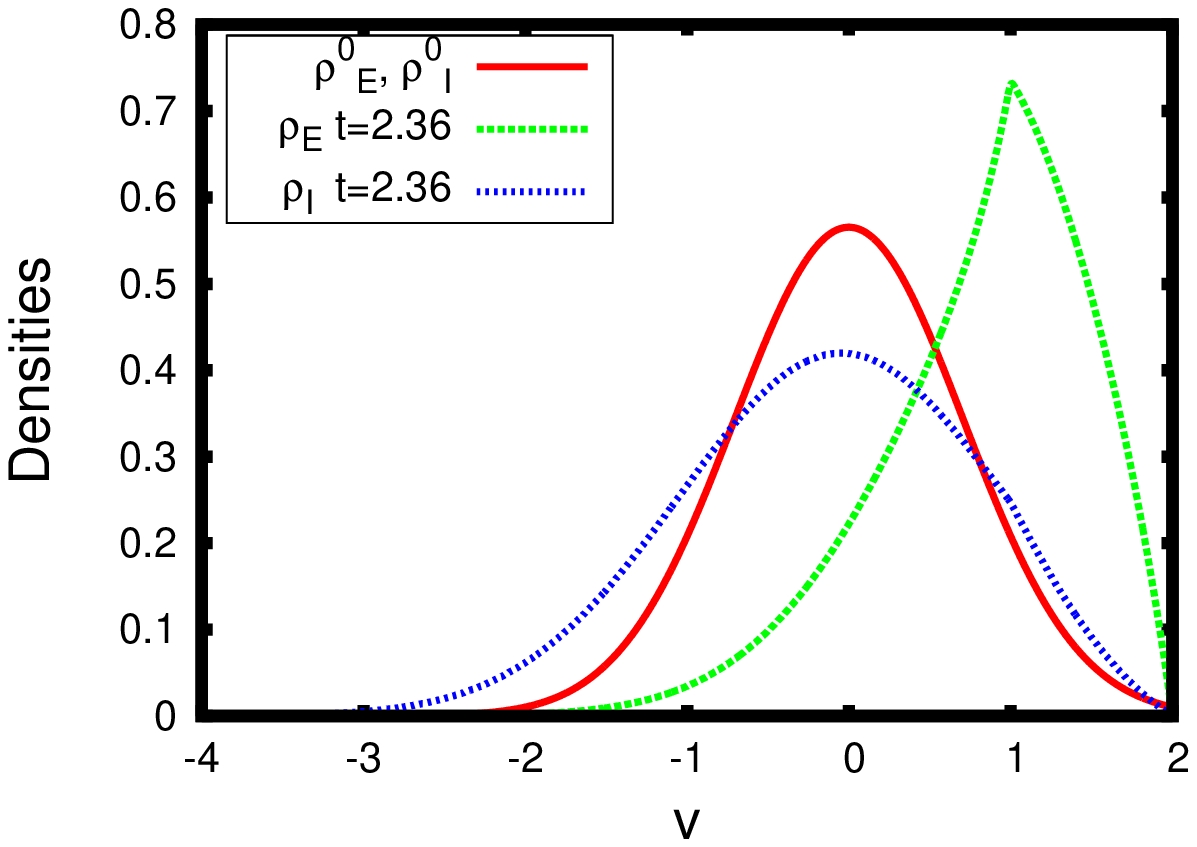}
\end{center}
\end{minipage}
\caption{Firing rates and probability densities for $b_E^E=3$, 
$b_I^E=0.75$, $b_E ^I=0.5$, $b_I^I=3$,  in case of a normalized Maxwellian 
initial condition with mean  0 and variance 0.5 (see \eqref{ci_maxwel}).The blow-up of $N_E$ cannot be avoided by
a large value of $b_I^I$.}
\label{blowup_bEE_bIE}
\end{figure}

For a purely inhibitory network the global
existence of its solutions was proved in \cite{CGGS}. 
Therefore, one could think that
 high values of $b_I^I$ could prevent the blow-up of the
excitatory-inhibitory system. However, Theorem \ref{th_blowup}
shows that this is not the case and 
 Fig. \ref{blowup_bEE_bIE} describes this situation;
although the value of  $b_I^I$ is big, a high value
of  $b_E^E$ causes the divergence of  $N_E(t)$
and the blow-up of the system.

\

\emph{\textbf{Steady states.}} 
The proof of Theorem \ref{th_eq} provides a strategy 
to find numerically the stationary firing rates, which
consists of finding the intersection points between
the functions $\mathcal{F}(N_E)$ and
constant 1 (see \eqref{ecF}). With this idea,  
 we have plotted  both functions 
for different parameter values.

\begin{figure}[H]
\begin{center}
\begin{minipage}[c]{0.33\linewidth}
\begin{center}
\includegraphics[width=\textwidth]{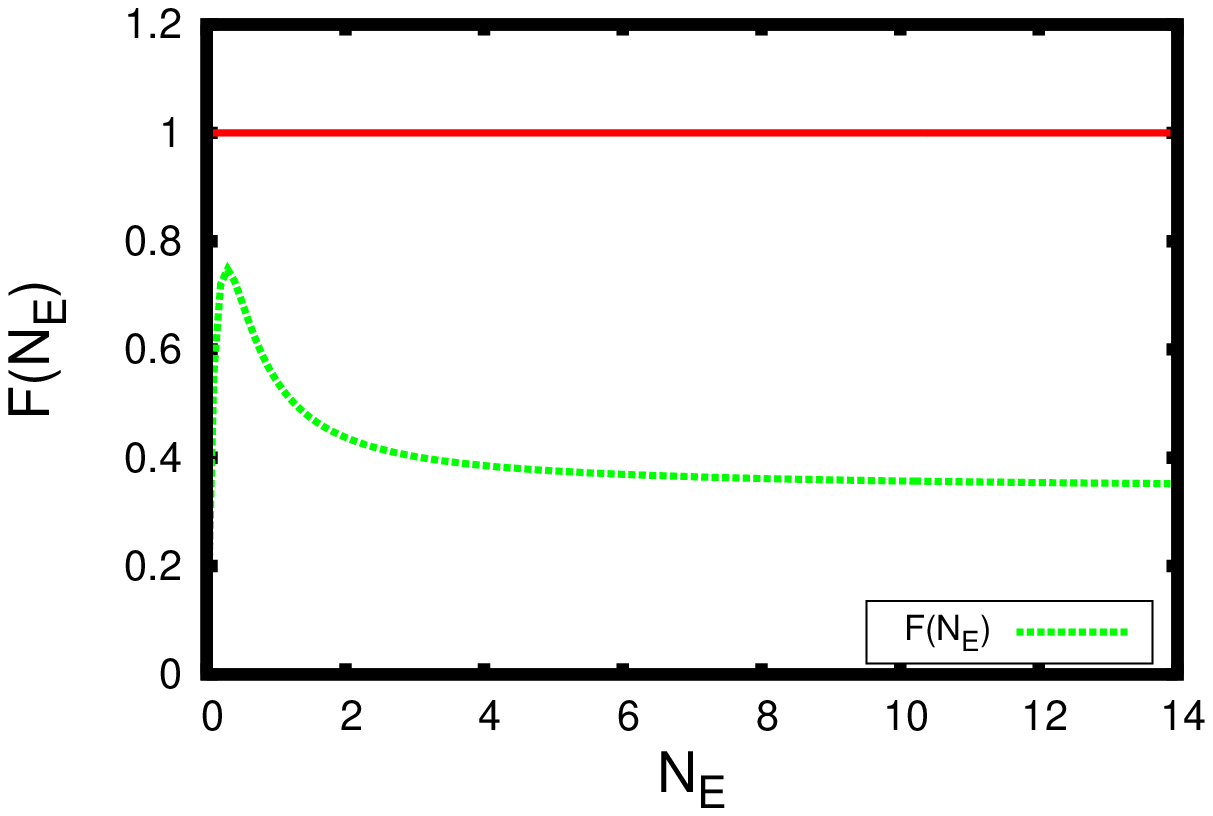}
\end{center}
\end{minipage}
\begin{minipage}[c]{0.33\linewidth}
\begin{center}
\includegraphics[width=\textwidth]{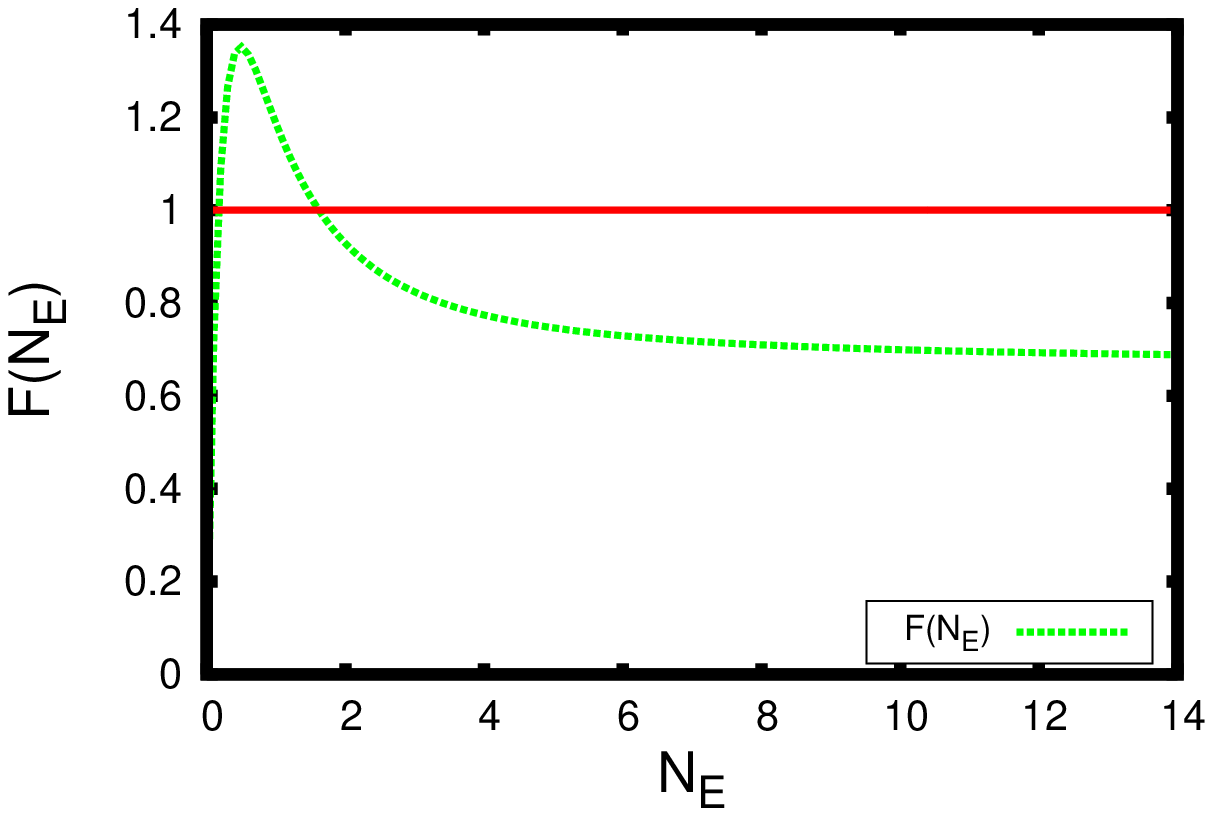}
\end{center}
\end{minipage}
\end{center}
\caption{$\mathcal{F}(N_E)$ for different parameter values 
corresponding to the first case of Theorem \ref{th_eq}:
 there are no steady states (left) or there is an even number of steady states (right).
\newline
Left figure:  
$b_E^E=3$, $b_I^E=0.75$, $b_E^I=0.5$ and $b_I^I=5$. 
Right figure:  $b_E^E=1.8$, $b_I^E=0.75$, $b_E^I=0.5$ and $b_I^I=0.25$.}
\label{caso1}
\end{figure}

The first case of Theorem \ref{th_eq} 
(there is an even number of steady states or there are no steady states) 
is shown in Fig. \ref{caso1}.
In this situation, the relation between the parameters implies 
$\displaystyle \lim_{N_E \to \infty} \mathcal{F}(N_E)<1$. 
In the left plot, the pure connectivity parameters
($b_E^E$ and $b_I^I$) are high in comparison with the connectivity
parameters $b_E^I$ and $b_I^E$, in such a way that there are no
steady states. 
In
 the right plot,
 there are two steady states because the pure
connectivity parameters are small, since in this case 
the maximum value of $\mathcal{F}$ is bigger than one.

\begin{figure}[H]
\begin{minipage}[c]{0.33\linewidth}
\begin{center}
\includegraphics[width=\textwidth]{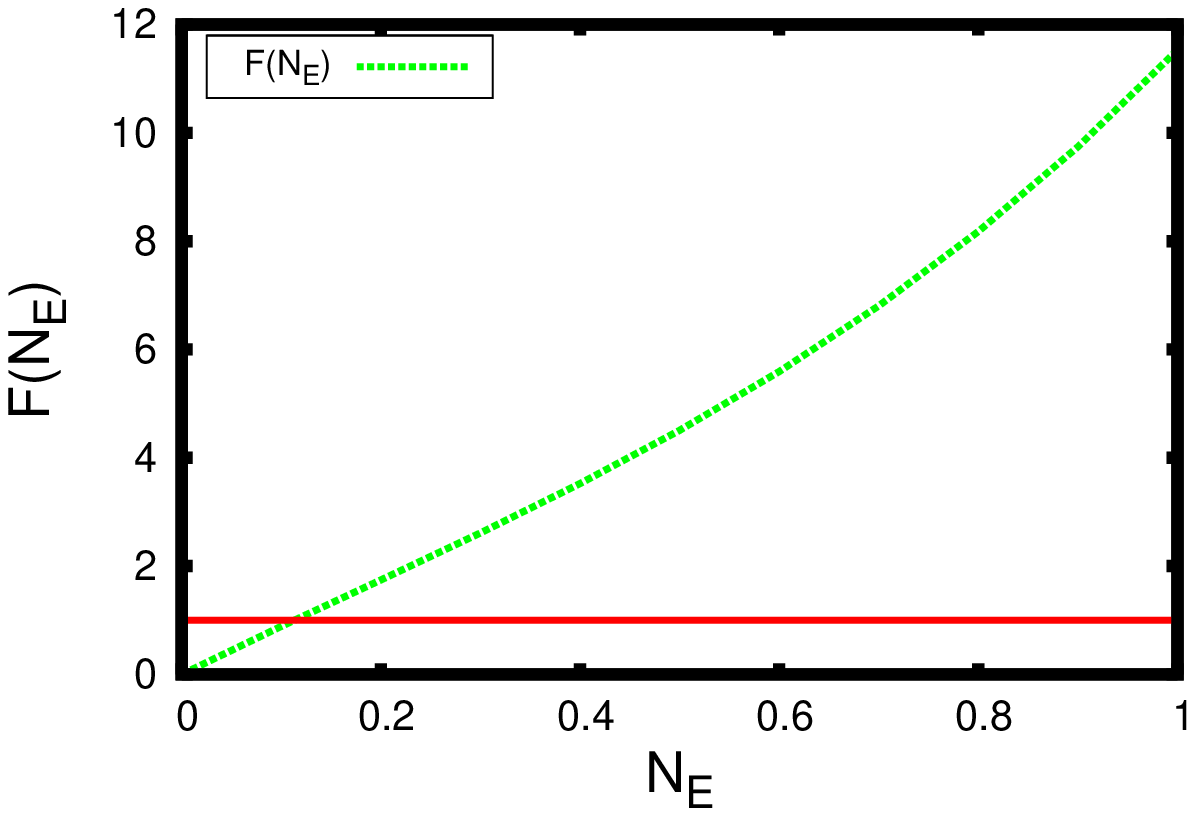}
\end{center}
\end{minipage}
\begin{minipage}[c]{0.33\linewidth}
\begin{center}
\includegraphics[width=\textwidth]{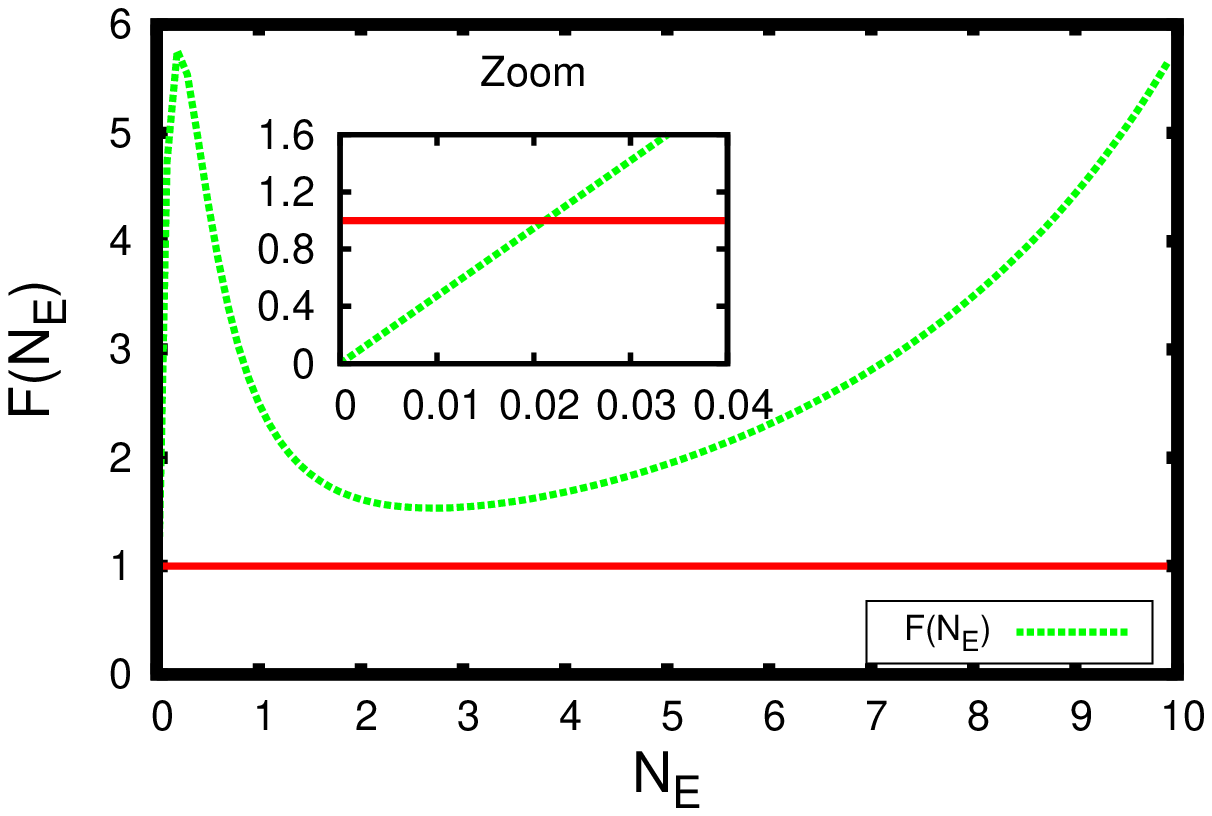}
\end{center}
\end{minipage}
\begin{minipage}[c]{0.33\linewidth}
\begin{center}
\includegraphics[width=\textwidth]{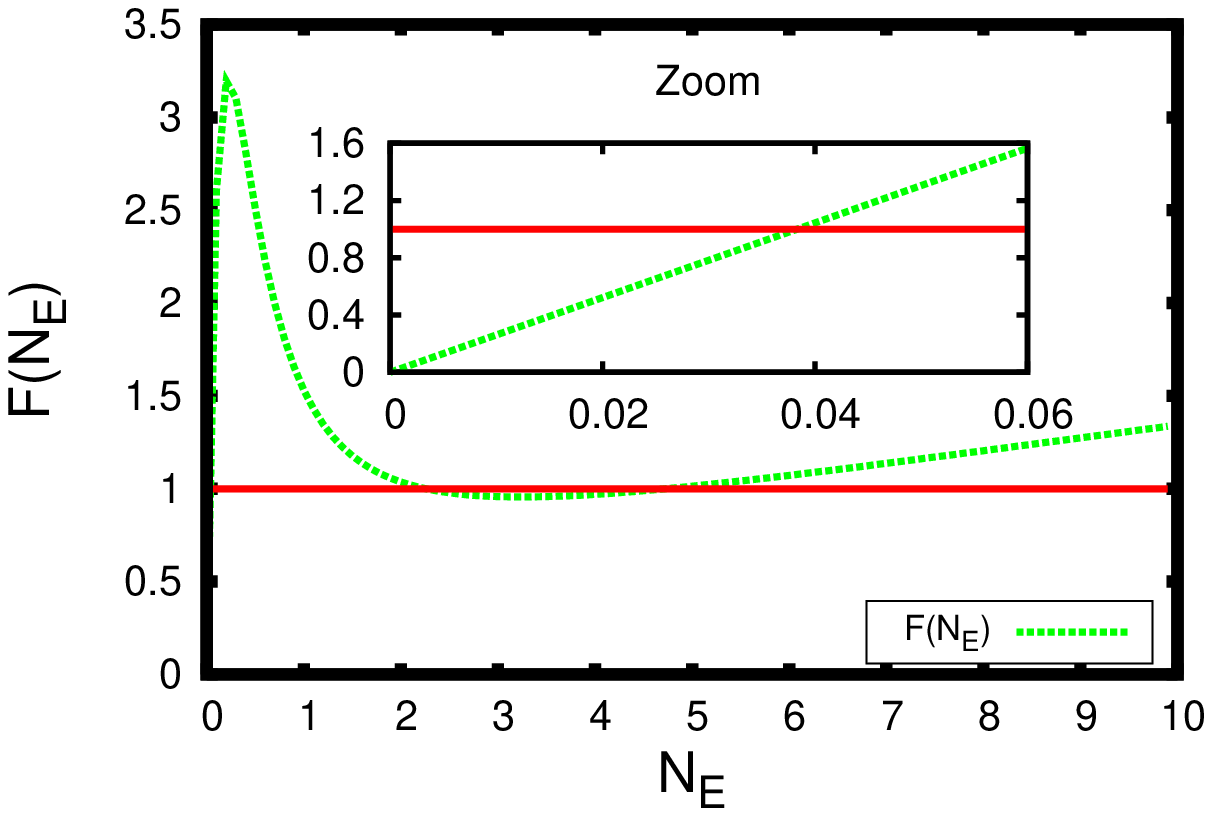}
\end{center}
\end{minipage}
\caption{$\mathcal{F}(N_E)$ for different parameter values 
corresponding to the second case of Theorem \ref{th_eq}:
there is an odd number of steady states.
\newline
Left figure:  $b_E^E=0.5$, $b_I^E=0.5$, $b_E^I=3$ and $b_I^I=0.5$ (one steady state).   
Center figure:  $b_E^E=3$,  $b_I^E=9$, $b_E^I=0.5$, $b_I^I=0.25$ (one steady state). 
Right figure: $b_E^E=3$,  $b_I^E=7$, $b_E^I=0.5$ and $b_I^I=0.25$ (three steady states). 
}
\label{caso2}
\end{figure}

The second case of Theorem \ref{th_eq} (there is an odd number of steady 
states), which implies  
$\displaystyle \lim_{N_E \to \infty} \mathcal{F}(N_E)>1$,  
is depicted in Fig. \ref{caso2}.
In the left plot, $b_E^E$ is small enough such that $\mathcal{F}$ is
an increasing function,  and therefore there is a unique steady state.
Also, in the center plot there is only one stationary solution, but
in this case the values of connectivity parameters do not guarantee 
the monotonicity of  $\mathcal{F}$.
Finally, the right plot shows values of connectivity parameters for which
there are three steady states.

\begin{figure}[H]
\begin{center}
\begin{minipage}[c]{0.33\linewidth}
\begin{center}
\includegraphics[width=\textwidth]{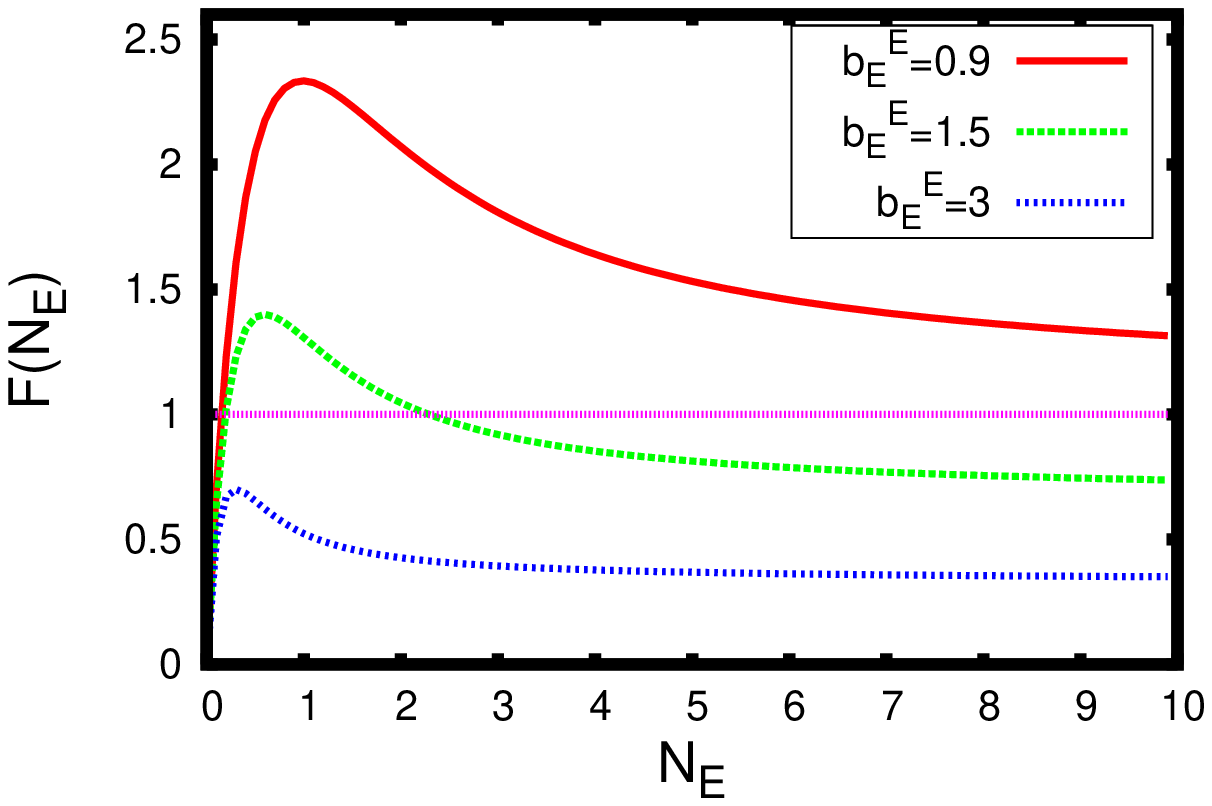}
\end{center}
\end{minipage}
\begin{minipage}[c]{0.33\linewidth}
\begin{center}
\includegraphics[width=\textwidth]{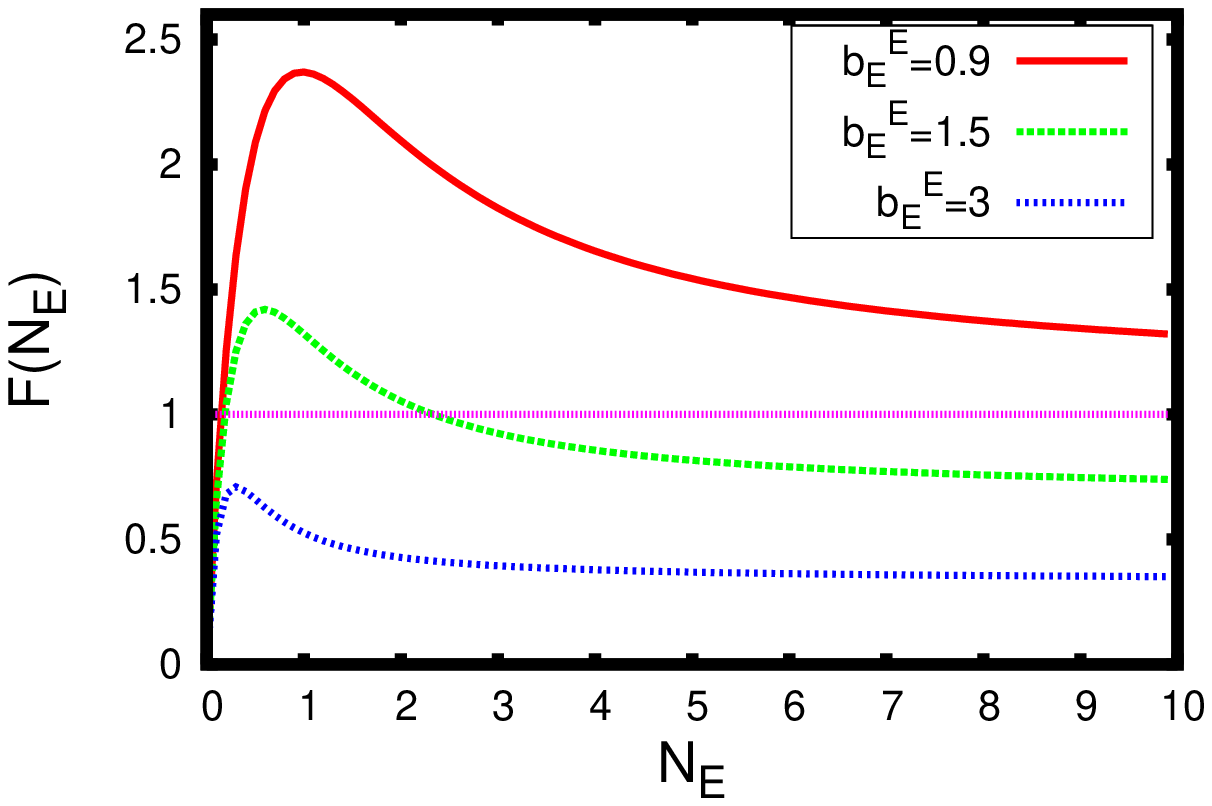}
\end{center}
\end{minipage}
\end{center}
\caption{Comparison between an uncoupled excitatory-inhibitory
network  ($b_E^I=b_I^E=0$) and a coupled network with small $b_E^I$ and 
$b_I^E$.  The qualitative behavior is the same in both cases.
\newline
Left figure:  $b_I^E=b_E^I=0$, $b_I^I=0.25$,  
and  different values for $b_E^E$.
Right figure:  $b_I^E=b_E^I=0.1$, $b_I^I=0.25$ and  
different values for $b_E^E$.}
\label{cruzados_chicos}
\end{figure}

\begin{figure}[H]
\begin{center}
\begin{minipage}[c]{0.33\linewidth}
\begin{center}
\includegraphics[width=\textwidth]{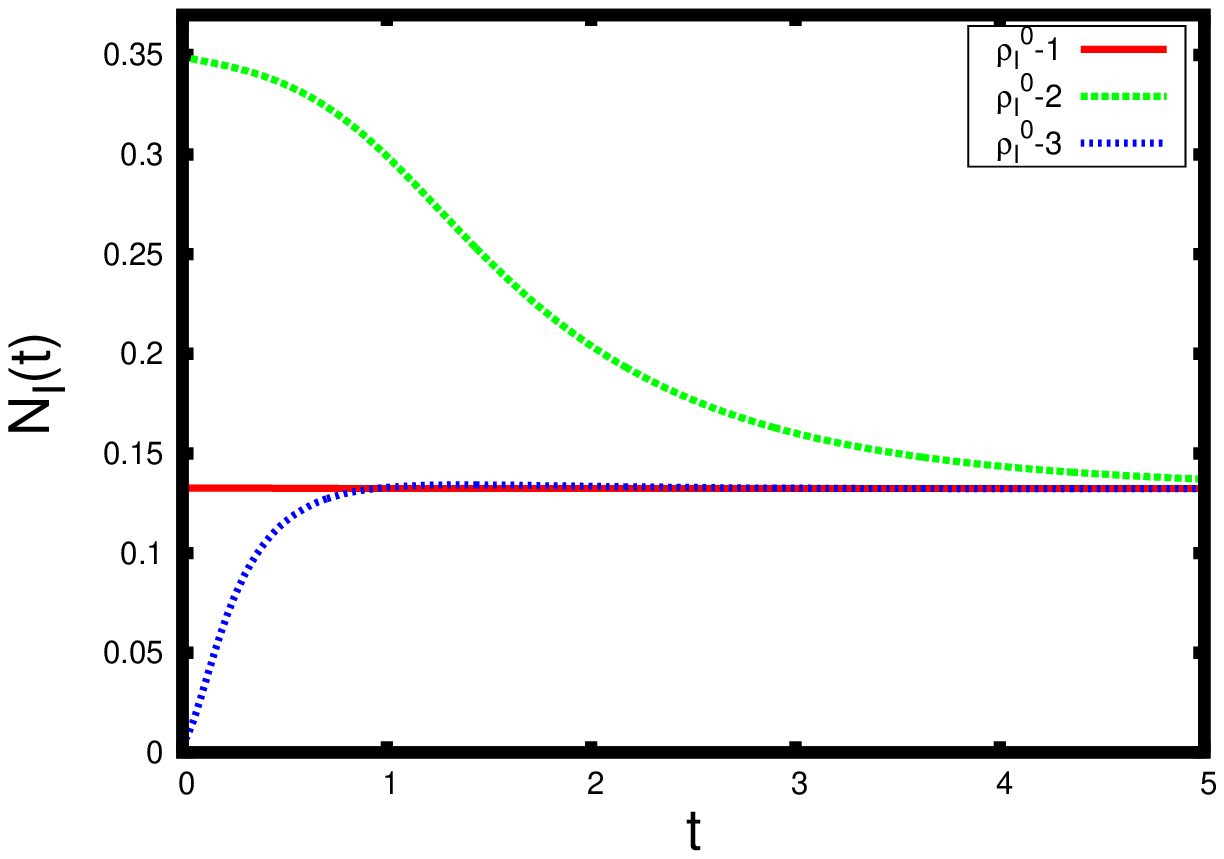}
\end{center}
\end{minipage}
\begin{minipage}[c]{0.33\linewidth}
\begin{center}
\includegraphics[width=\textwidth]{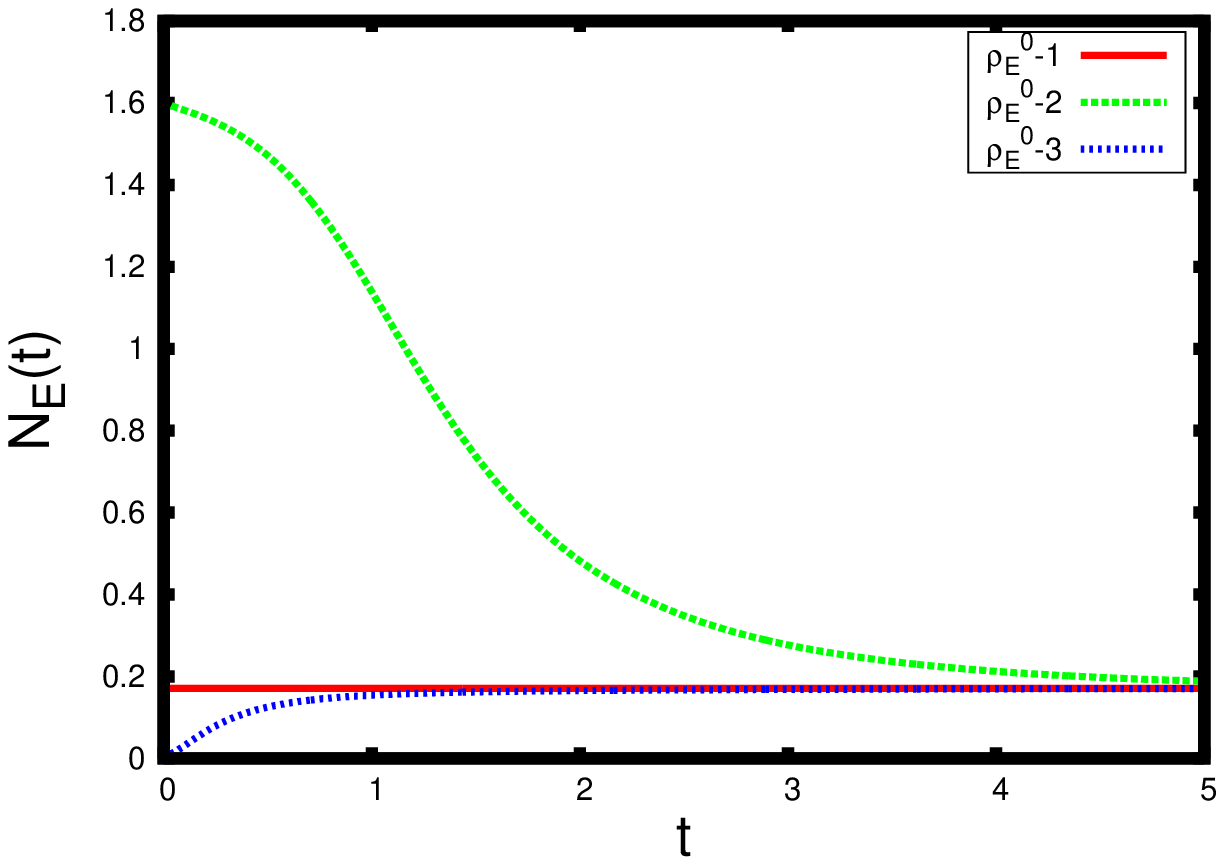}
\end{center}
\end{minipage}
\end{center}
\caption{Firing rates for the case of two steady states 
(right in Fig. \ref{caso1}), for different initial conditions: 
$\rho_{\alpha}^0-1,2$ are given by the profile 
\eqref{soleq} with $(N_E,N_I)$  stationary values and 
$\rho_{\alpha}^0-3$ is a normalized Maxwellian with mean 0 and variance 0.25
(see \eqref{ci_maxwel}). For both firing rates, the lower steady state seems to be asymptotically stable  whereas the higher one seems to be unstable.}
\label{estabilidad2}
\end{figure}

To conclude the analysis of the number of steady states, 
in Fig. \ref{cruzados_chicos} it is depicted a comparison between 
an uncoupled excitatory-inhibitory
network  ($b_E^I=b_I^E=0$) and a coupled network with small $b_E^I$ and $b_I^E$.
In Fig. \ref{cruzados_chicos} (left) the number of steady states for 
the uncoupled system ($b_E^I=b_I^E=0$) is analyzed, 
while Fig. \ref{cruzados_chicos} (right) investigates the number of steady 
states when the parameters $b_E^I$ and $b_I^E$  are small ($b_E^I=b_I^E=0.1$).
As expected, it can be observed that the number of steady states 
does not change if we choose small values for $b_E^I$ and $b_I^E$.
It depends only on the value of $b_E^E$ in the same 
way as described in  \cite{CCP}. For small values of
$b_E^E$ there is a unique steady state. As $b_E^E$ increases,
 it appears another steady state, that merges with the first one and
then disappears, apparently yielding a saddle node bifurcation. 
In Fig. \ref{estabilidad2} we show the time evolution  of the firing
rates, for the case described in Fig. \ref{caso1} (right),
with two steady states. We use as initial condition profiles 
like the one presented in \eqref{soleq}.
The numerical results show that the  larger steady state is unstable,
while the lower one seems to be stable.
Therefore, numerical stability analysis 
indicates that the unique steady state is stable,
while, when there are two steady states the highest seems to be unstable.

\begin{figure}[H]
\begin{center}
\begin{minipage}[c]{0.7\linewidth}
\begin{center}
\includegraphics[width=\textwidth]{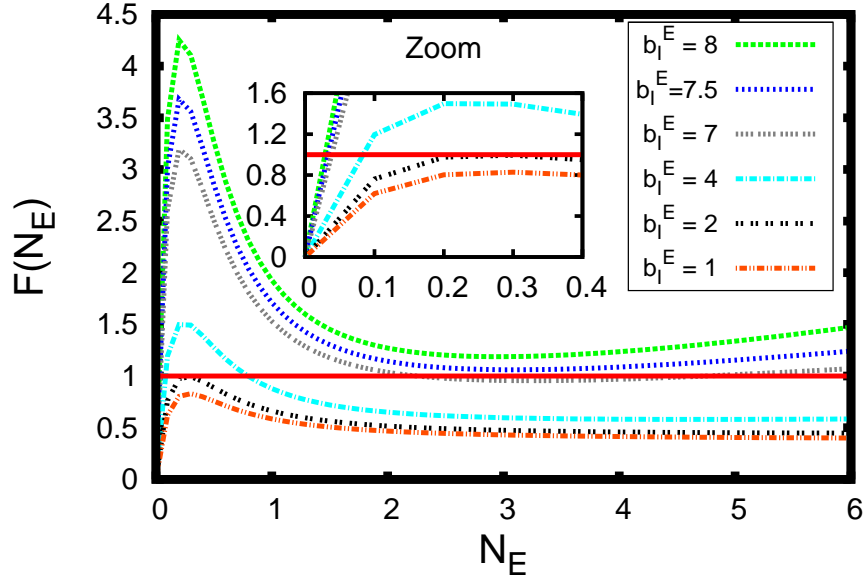}
\end{center}
\end{minipage}
\end{center}
\caption{Analysis of the number of steady states for $b_E^E=3$, $b_E^I=0.5$, $b_I^I=0.25$ and different values for $b_I^E$.}
\label{bifurcacion}
\end{figure}

\begin{figure}[H]
\begin{minipage}[c]{0.33\linewidth}
\begin{center}
\includegraphics[width=\textwidth]{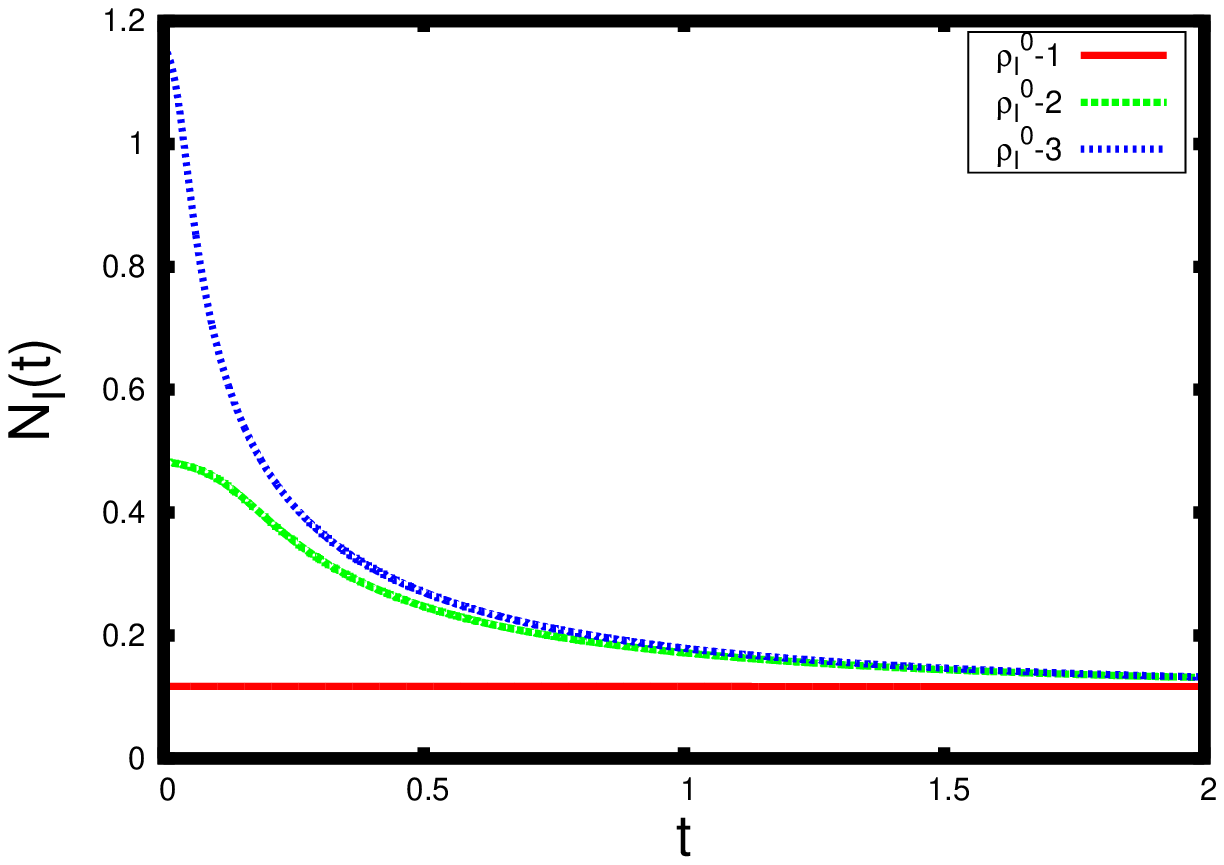}
\end{center}
\end{minipage}
\begin{minipage}[c]{0.33\linewidth}
\begin{center}
\includegraphics[width=\textwidth]{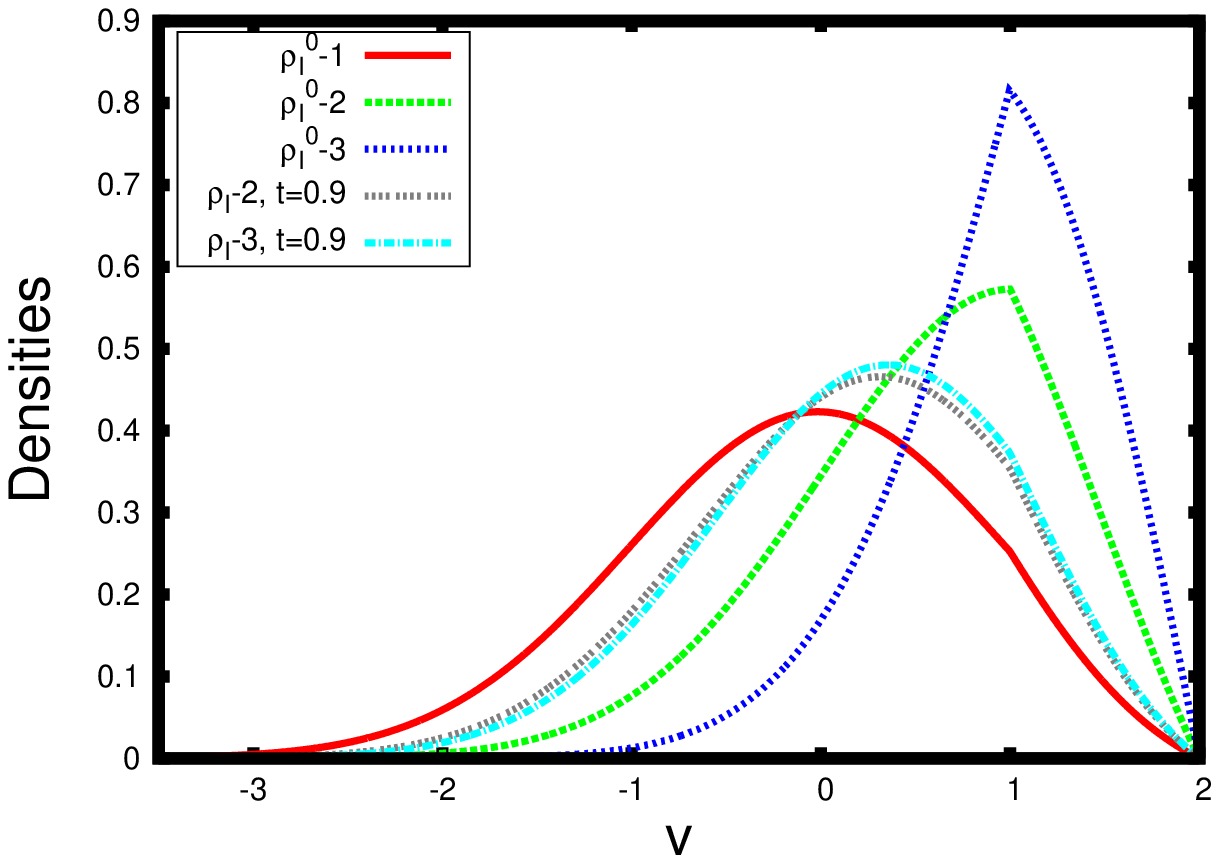}
\end{center}
\end{minipage}
\begin{minipage}[c]{0.33\linewidth}
\begin{center}
\includegraphics[width=\textwidth]{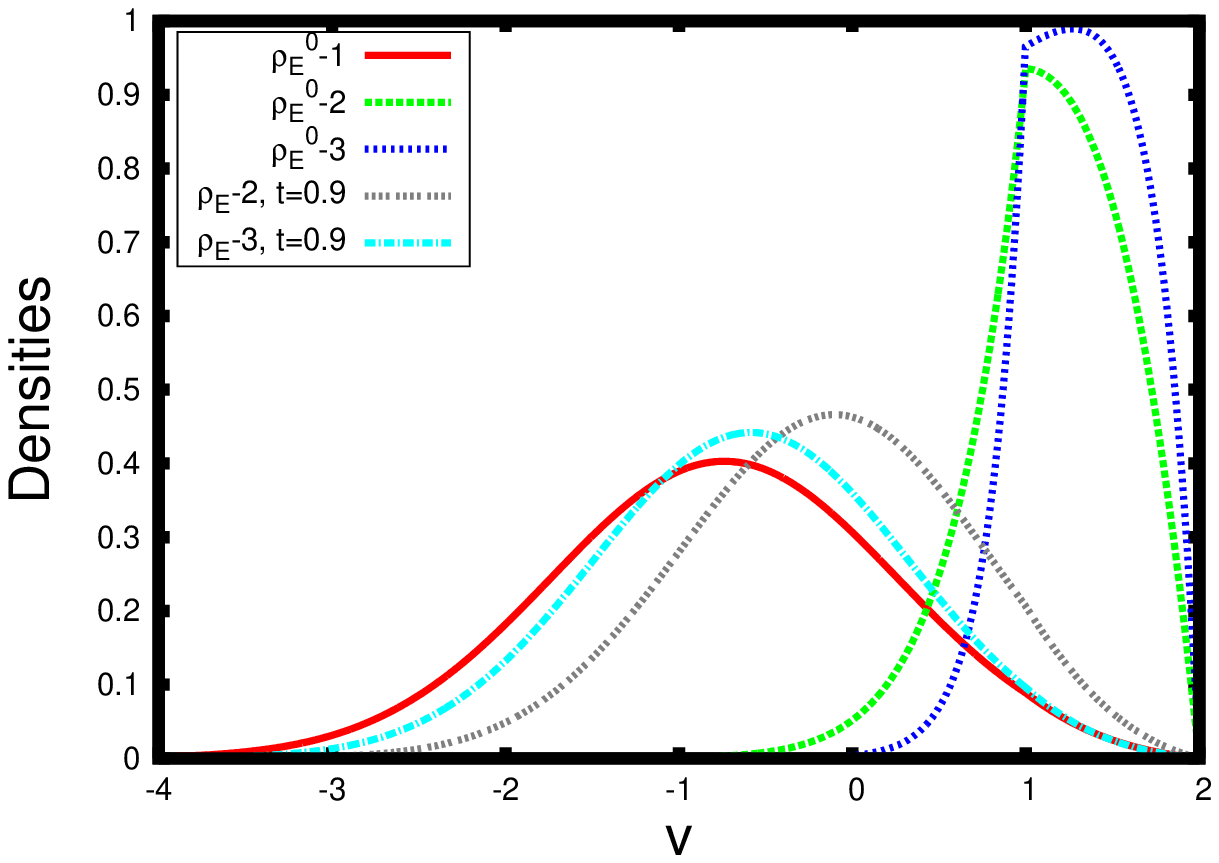}
\end{center}
\end{minipage}
\caption{Stability analysis for the case of three steady states 
(right in Fig. \ref{caso2}). 
\newline
Left figure: firing rates for different initial conditions: 
$\rho_{\alpha}^0-1,2,3$ which are given by the profile 
\eqref{soleq} with $(N_E,N_I)$ stationary values.Only the lowest steady state seems to be asymptotically stable. 
Center and right figure: evolution of the probability densities.
\newline
(Simulations were developed considering $v\in [-6,2]$).}
\label{estabilidad3}
\end{figure}

In this direction, another interesting bifurcation 
analysis that can be obtained in terms of the parameter $b_I^E$ is 
depicted in Fig. \ref{bifurcacion}. For large values of $b_E^I$ 
there is a unique stable steady. Then, as it decreases it appears 
another  steady state that bifurcates and gives rise to two equilibria. 
The largest  disappears, and the lower one approaches the smallest
 steady state. Afterwards, both  dis\-appear, probably through a
 saddle node bifurcation. Numerical stability analysis determines 
that the lowest steady state is always stable, while the other ones 
are all unstable.
The Fig. \ref{caso2} (right) depicts three steady states.
The stability analysis of this situation is more complicated than
the previous one (the case of two stationary solutions).
Fig.  \ref{estabilidad3} shows the time evolution of the solutions
for different initial data; the steady state with less firing rate
seems to be stable, while the other two steady states are unstable.

\

In conclusion, in the present work we have extended the known results for purely excitatory or inhibitory networks \cite{CCP} to excitatory-inhibitory coupled networks. We have proved that the presence of an inhibitory population in a coupled network does not avoid the blow-up phenomenon, as happens in a purely inhibitory network. Besides, we have analysed the number of steady states of the system, which is a more complicate issue than in the case of uncoupled systems. For small connectivity pa\-ra\-me\-ter values, we have shown that solutions converge exponentially fast to the unique steady state. Our analytical and numerical results contribute to support the NNLIF system as an appropriate model to describe well known neurophysiological phenomena, as for example synchronization/asynchronization of the network, since the blow-up in finite time might depict a synchronization of a part of the network, while the presence of a unique asymptotically stable stationary solution represents an asynchronization of the network. In addition, the abundance in the number of steady states, in terms of the 
connectivity parameter values, 
that can be observed for this simplified model,  
probably will help us to characterize situations of multi-stability 
for more complete NNLIF models and also other
models including conductance variables as in \cite{CCTa}. 
In \cite{CP} it was shown that if  a
 refractory period is included in the model, there are situations of
 multi-stability, with two stable and one unstable steady state.
In \cite{CCTa} bi-stability phenomena were  numerically described. 
Multi-stable networks are
 related, for instance, to the visual 
perception and the decision making \cite{GrayAndSinger:89,AD09}.
In these directions, some tasks which remain open for 
the  NNLIF model are: The  analytical study of the 
stability of the network when there is more than one steady state,
the possible presence of periodic solutions and the existence of
solutions for coupled networks, maybe following a similiar
 strategy as in \cite{GG09}. 
A more extensive bifurcation analysis could probably  contribute to find answers to these questions.

\

{\footnotesize \noindent\textit{Acknowledgments.} The 
authors acknowledge support from the project MTM2014-52056-P 
of Spanish Ministerio de Econom\'\i a y Competitividad and 
the European Regional Development Fund (ERDF/FEDER).} 

\newpage

\bibliographystyle{siam}
\bibliography{Bib.bib}

\begin{thebibliography}{10}

\bibitem{AD09}
{\sc L.~Albantakis and G.~Deco}, {\em The encoding of alternatives in
  multiple-choice decision making}, Proc Natl Acad Sci U S A, 106 (2009),
  pp.~10308--10313.

\bibitem{BrGe}
{\sc R.~Brette and W.~Gerstner}, {\em Adaptive exponential integrate-and-fire
  model as an effective description of neural activity}, Journal of
  neurophysiology, 94 (2005), pp.~3637--3642.

\bibitem{brunel}
{\sc N.~Brunel}, {\em Dynamics of sparsely connected networks of excitatory and
  inhibitory spiking networks}, J. Comp. Neurosci., 8 (2000), pp.~183--208.

\bibitem{BrHa}
{\sc N.~Brunel and V.~Hakim}, {\em Fast global oscillations in networks of
  integrate-and-fire neurons with long firing rates}, Neural Computation, 11
  (1999), pp.~1621--1671.

\bibitem{CCP}
{\sc M.~J. C\'aceres, J.~A. Carrillo, and B.~Perthame}, {\em Analysis of
  nonlinear noisy integrate $\&$ fire neuron models: blow-up and steady
  states}, Journal of Mathematical Neuroscience, 1-7 (2011).

\bibitem{CCTa}
{\sc M.~J. C\'aceres, J.~A. Carrillo, and L.~Tao}, {\em A numerical solver for
  a nonlinear fokker-planck equation representation of neuronal network
  dynamics}, J. Comp. Phys., 230 (2011), pp.~1084--1099.

\bibitem{CP}
{\sc M.~J. C\'aceres and B.~Perthame}, {\em Beyond blow-up in excitatory
  integrate and fire neuronal networks: refractory period and spontaneous
  activity}, Journal of theoretical Biology, 350 (2014), pp.~81--89.

\bibitem{carrillo2014qualitative}
{\sc J.~Carrillo, B.~Perthame, D.~Salort, and D.~Smets}, {\em Qualitative
  properties of solutions for the noisy integrate \& fire model in
  computational neuroscience}, Nonlinearity, 25 (2015), pp.~3365--3388.

\bibitem{CGGS}
{\sc J.~A. Carrillo, M.~d.~M. Gonz\'alez, M.~P. Gualdani, and M.~E. Schonbek},
  {\em Classical solutions for a nonlinear fokker-planck equation arising in
  computational neuroscience}, Comm. in Partial Differential Equations, 38
  (2013), pp.~385--409.

\bibitem{chevallier2015mean}
{\sc J.~Chevallier}, {\em Mean-field limit of generalized hawkes processes},
  arXiv preprint arXiv:1510.05620,  (2015).

\bibitem{chevallier2015microscopic}
{\sc J.~Chevallier, M.~J. C{\'a}ceres, M.~Doumic, and P.~Reynaud-Bouret}, {\em
  Microscopic approach of a time elapsed neural model}, Mathematical Models and
  Methods in Applied Sciences, 25 (2015), pp.~2669--2719.

\bibitem{delarue2015particle}
{\sc F.~Delarue, J.~Inglis, S.~Rubenthaler, and E.~Tanr{\'e}}, {\em Particle
  systems with a singular mean-field self-excitation. application to neuronal
  networks}, Stochastic Processes and their Applications, 125 (2015),
  pp.~2451--2492.

\bibitem{delarue2015global}
{\sc F.~Delarue, J.~Inglis, S.~Rubenthaler, E.~Tanr{\'e}, et~al.}, {\em Global
  solvability of a networked integrate-and-fire model of mckean--vlasov type},
  The Annals of Applied Probability, 25 (2015), pp.~2096--2133.

\bibitem{Henry:13}
{\sc G.~Dumont and J.~Henry}, {\em Synchronization of an excitatory
  integrate-and-fire neural network}, Bull. Math. Biol., 75 (2013),
  pp.~629--648.

\bibitem{dumont2015noisy}
{\sc G.~Dumont, J.~Henry, and C.~O. Tarniceriu}, {\em Noisy threshold in
  neuronal models: connections with the noisy leaky integrate-and-fire model},
  arXiv preprint arXiv:1512.03785,  (2015).

\bibitem{GK}
{\sc W.~Gerstner and W.~Kistler}, {\em Spiking neuron models}, Cambridge Univ.
  Press, Cambridge, 2002.

\bibitem{GG09}
{\sc M.~d.~M. Gonz\'alez and M.~P. Gualdani}, {\em Asymptotics for a symmetric
  equation in price formation}, App. Math. Optim., 59 (2009), pp.~233--246.

\bibitem{GrayAndSinger:89}
{\sc C.~M. Gray and W.~Singer}, {\em {Stimulus-specific neuronal oscillations
  in orientation columns of cat visual cortex}}, Proc Natl Acad Sci U S A, 86
  (1989), pp.~1698--1702.

\bibitem{G}
{\sc T.~Guillamon}, {\em An introduction to the mathematics of neural
  activity}, Butl. Soc. Catalana Mat., 19 (2004), pp.~25--45.

\bibitem{mg}
{\sc M.~Mattia and P.~Del~Giudice}, {\em Population dynamics of interacting
  spiking neurons}, Phys. Rev. E, 66 (2002), p.~051917.

\bibitem{omurtag}
{\sc A.~Omurtag, K.~B. W., and L.~Sirovich}, {\em On the simulation of large
  populations of neurons}, J. Comp. Neurosci., 8 (2000), pp.~51--63.

\bibitem{PPD}
{\sc K.~Pakdaman, B.~Perthame, and D.~Salort}, {\em Dynamics of a structured
  neuron population}, Nonlinearity, 23 (2010), pp.~55--75.

\bibitem{PPD2}
{\sc K.~Pakdaman, B.~Perthame, and D.~Salort}, {\em Relaxation and
  self-sustained oscillations in the time elapsed neuron network model}, SIAM
  Journal on Applied Mathematics, 73 (2013), pp.~1260--1279.

\bibitem{pakdaman2014adaptation}
{\sc K.~Pakdaman, B.~Perthame, and D.~Salort}, {\em Adaptation and fatigue
  model for neuron networks and large time asymptotics in a nonlinear
  fragmentation equation}, The Journal of Mathematical Neuroscience (JMN), 4
  (2014), pp.~1--26.

\bibitem{perthame2013voltage}
{\sc B.~Perthame and D.~Salort}, {\em On a voltage-conductance kinetic system
  for integrate and fire neural networks}, Kinetic and related models, AIMS, 6
  (2013), pp.~841--864.

\bibitem{Caikinetic}
{\sc A.~V. Rangan, G.~Kova{\u c}i{\u c}, and D.~Cai}, {\em Kinetic theory for
  neuronal networks with fast and slow excitatory conductances driven by the
  same spike train}, Physical Review E, 77 (2008), pp.~1--13.

\bibitem{RBW}
{\sc A.~Renart, N.~Brunel, and X.-J. Wang}, {\em Mean-field theory of
  irregularly spiking neuronal populations and working memory in recurrent
  cortical networks}, in Computational Neuroscience: A comprehensive approach,
  J.~Feng, ed., Chapman \& Hall/CRC Mathematical Biology and Medicine Series,
  2004.

\bibitem{Risken}
{\sc H.~Risken}, {\em The Fokker-Planck Equation: Methods of solution and
  approximations}, 2nd. edn. Springer Series in Synergetics, vol 18.
  Springer-Verlag, Berlin, 1989.

\bibitem{RB}
{\sc C.~Rossant, D.~F.~M. Goodman, B.~Fontaine, J.~Platkiewicz, A.~K.
  Magnusson, and R.~Brette}, {\em Fitting neuron models to spike trains},
  Frontiers in Neuroscience, 5 (2011), pp.~1--8.

\bibitem{Touboul_2008}
{\sc J.~Touboul}, {\em Bifurcation analysis of a general class of nonlinear
  integrate-and-fire neurons}, SIAM J. Appl. Math., 68 (2008), pp.~1045--1079.

\bibitem{Touboul_AQIF}
{\sc J.~Touboul}, {\em Importance of the cutoff value in the quadratic adaptive
  integrate-and-fire model}, Neural Computation, 21 (2009), pp.~2114--2122.

\bibitem{T}
{\sc H.~Tuckwell}, {\em Introduction to Theoretical Neurobiology}, Cambridge
  Univ. Press, Cambridge, 1988.

\end{thebibliography}

\emph{E-mail address:} \texttt{caceresg@ugr.es}

\emph{E-mail address:} \texttt{ricardaschneider@ugr.es}

\end{document}